\documentclass[preprint]{elsarticle}

\usepackage{graphicx}
\usepackage{lineno}
\usepackage{amssymb}
\usepackage{mathrsfs}
\usepackage{mhchem}
\usepackage{subcaption}
\usepackage{epstopdf}
\usepackage{float}
\usepackage{bbold}
\usepackage{amsmath}
\usepackage{color} 
\usepackage{ulem} 
\usepackage{algorithmicx}
\usepackage{algorithm}
\usepackage{algpseudocode}

\usepackage{caption}
\usepackage{hyperref}
\usepackage{amsfonts}
\usepackage{graphicx}

\newcommand{\R}{\mathbb{R}}

\renewcommand{\u}{\mathbf{u}}
\newcommand{\U}{\mathbf{U}}
\renewcommand{\v}{\mathbf{v}}

\renewcommand{\d}{\mathbf{d}}
\newcommand{\g}{\mathbf{g}}
\newcommand{\f}{\mathbf{f}}

\newcommand{\F}{\mathbf{F}}

\renewcommand{\c}{\mathbf{c}}

\numberwithin{theorem}{section}
\usepackage{amsopn}

\date{\today}

\begin{document}

\begin{frontmatter}
\title{A Study of the Numerical Stability of an ImEx Scheme with Application to the Poisson-Nernst-Planck Equations }

\author{M.C. Pugh\corref{mycorrespondingauthor}}
\address{Department of Mathematics, University of Toronto, \\
40 St George St, Toronto, ON  M5S 2E4, Canada}
\cortext[mycorrespondingauthor]{Corresponding author}
\ead{mpugh@math.utoronto.ca}

\author{David Yan}
\address{Department of Electrical and Computer Engineering, University of Toronto}

\author{F.P. Dawson}
\address{Department of Electrical and Computer Engineering, University of Toronto}


\begin{abstract}
The Poisson-Nernst-Planck equations with generalized
  Frumkin-Butler-Volmer boundary conditions (PNP-FBV) describe ion
  transport with Faradaic reactions and have applications in a wide
  variety of fields.
 We solve the PNP-FBV equations 
using an adaptive time-stepper based
  on a second-order variable step-size, semi-implicit, backward
differentiation formula (VSSBDF2).  
When the
  underlying dynamics are such that the solutions converge to
  a steady-state solution, we observe that the adaptive time-stepper produces
  solutions that ``nearly'' converge to the steady state and that,
  simultaneously, the time-step sizes stabilize to a limiting size
  $dt_\infty$.
  Linearizing the SBDF2 scheme about the steady
  state solution, we demonstrate that the linearized scheme is conditionally
  stable and that this is the cause of the adaptive time-stepper's
  behaviour.   Mesh-refinement, as well as a study of the eigenvectors
  corresponding to the critical eigenvalues, demonstrate that the conditional stability is not due to a time-step
  restriction caused by high-frequency contributions. We study the stability domain of the linearized scheme and
find that its boundary can have corners as well as jump discontinuities.  
A jump discontinuity means there can be parameter values,
$\epsilon_1$ and $\epsilon_2$, that are very close to one another and a
time-step size $dt$ so that the computation of the $\epsilon_1$ problem
is stable and the computation of the $\epsilon_2$ problem is unstable.
\end{abstract}



\end{frontmatter}


\section{Introduction}\label{introduction}


 The Poisson-Nernst-Planck
(PNP) equations are a parabolic-elliptic system that models the transport of charged species subject to
diffusion and electromigration.  The generalized
Frumkin-Butler-Volmer (FBV) boundary conditions are
nonlinear boundary conditions that model
chemical reactions at the electrodes.  
The PNP equations model the behaviour in the bulk; the 
electrodes are located at the boundary of the ``bulk''
domain.
There is a singular perturbation parameter $\epsilon$; small values of
$\epsilon$ lead to thin boundary layers with sharp transitions to the
behaviour in the bulk.

The PNP-FBV system is both nonlinear and diffusive.  Hence 
a semi-implicit (also known as implicit-explicit) time-stepping scheme is a natural approach 
to take in hopes of avoiding stability restrictions on the time-step size
while also avoiding the computational slowness caused by having to solve
nonlinear equations.  The linear diffusive term is handled implicitly and the nonlinear
terms are handled explicitly. 
In this article, we study the stability properties of a
second-order semi-implicit backwards differencing formula (SBDF2)
as applied to the PNP-FBV system and find that the SBDF2 scheme
becomes conditionally stable as the underlying solution of PNP-FBV system
equilibrates.  
%

The PNP-FBV system can be forced at one of the boundaries using either an imposed voltage or imposed current.  Simulations considered a variety of imposed voltages and imposed currents; the adaptive time-stepper was vital in that it could refine, and subsequently coarsen,
the time-steps in response to fast changes in the imposed forcing  \cite{YanThesis,Yan2017}.
The adaptive time-stepper 
presented in  \cite{YanThesis} is  
based on a
second-order variable step-size, semi-implicit, backward
differentiation formula (VSSBDF2 \cite{Wang2008}).  In the companion
article \cite{YPD_Part1}, the adaptive VSSBDF2 time-stepper is presented in 
full and its speed and stability properties are compared to those of 
an adaptive time-stepper based on a
second-order variable step-size, {\it fully-implicit}, backward
differentiation formula (VSBDF2  \cite{Wang2008}).
It is demonstrated that for ``large'' values of the singular perturbation
parameter $\epsilon$, the (semi-implicit) VSSBDF2 adaptive time-stepper is faster
and for ``small'' values of $\epsilon$ the (fully-implicit) VSBDF2 adaptive time-stepper is faster.

\begin{figure}[htb!]
\centering
\includegraphics[width=0.7\linewidth]{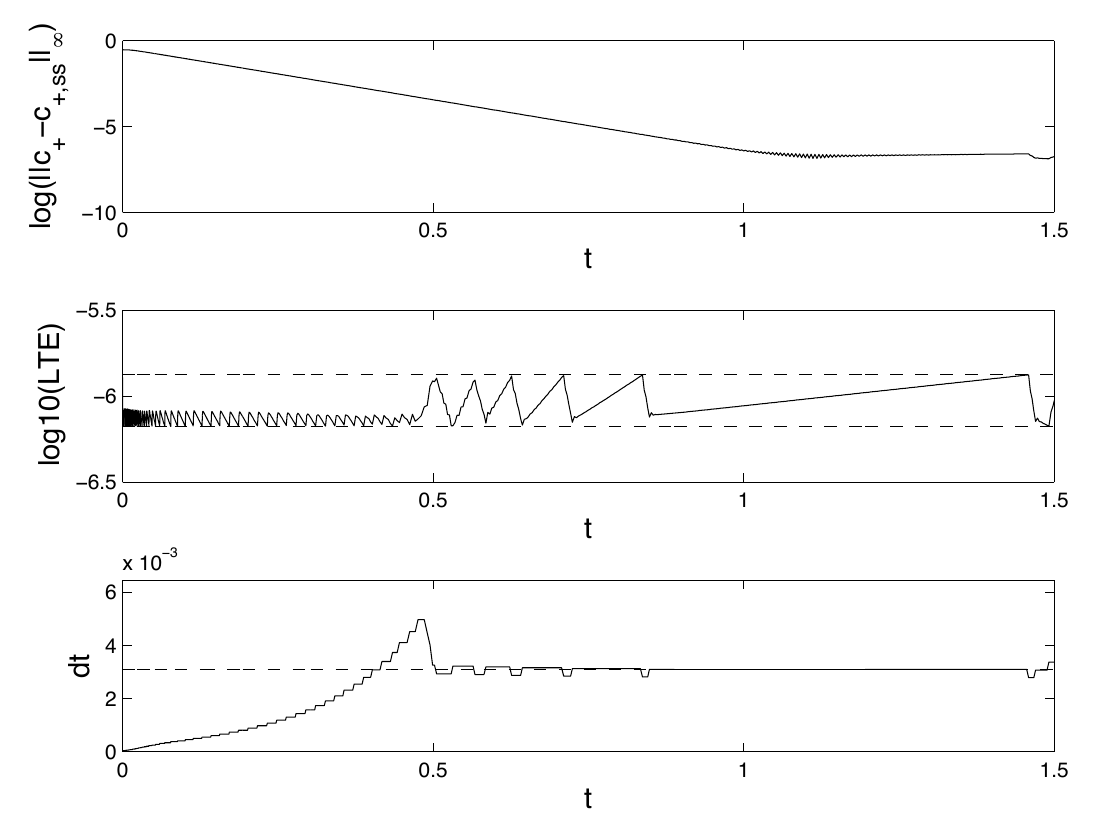}
\caption{PNP-FBV system \eqref{concentration_nondim}--\eqref{phi_bc_nondim_R} 
with constant imposed voltage $v(t) = 2$, $\epsilon = .05$, and all other physical parameters set to $1$.  The spatial mesh is uniform, $dx = 1/90$.
The initial data is $c_\pm(x,0) = 1 + .1 \, \sin(2 \pi x)$ and $\phi_x(1,0) = 0$.
The VSSBDF2 adaptive time-stepper is used.
\underline{Top plot:}  Comparison of the solutions found by the VSSBDF2 adaptive time-stepper to the pre-computed
numerical steady-state solution.  Plot is
$\log(\|\c_+^n - \c_{+,ss}\|_\infty)$ versus $t_n$.  Deviations of $\c_-^n$ and $\pmb{\phi}^n$ from the corresponding steady
state profiles behave similarly, as do their time derivatives as approximated
using \eqref{sbdf2}.  
\underline{Middle plot:} The logarithm of the
approximate local truncation error, \eqref{final_error}, is plotted versus time.
The dashed lines indicate the constraints set by the 
adaptive time-stepper:
$\log(tol \pm range)$.  See Appendix \ref{ATS}.
\underline{Bottom plot:} Time-step size, $dt$, plotted versus time.  The
timestep sizes stabilize to a limiting value, denoted $dt_\infty$.  The 
dashed line indicates the
stability restriction $dt^* = 3.1000 \times 10^{-3}$
computed using the linear stability analysis presented in Section \ref{linear_stability}. 
}
\label{Full_Model_tmax_effect}
\end{figure}

When the imposed forcing is
held constant, and the underlying physical solution relaxes to a steady-state solution, the expected behaviour of an adaptive time-stepper is that the
time-steps will grow until they reach the user-specified maximum time-step size, $dt_\text{max}$, and the numerical 
solution will converge to the
numerical steady state.   
This is precisely what is observed with the VSBDF2 adaptive time-stepper based on the
fully-implicit BDF2 scheme \cite{YPD_Part1}.
However, 
we found that 
the VSSBDF2 adaptive time-stepper based on the semi-implicit SBDF2 scheme did not
behave in the expected manner.  Instead, as shown in Figure \ref{Full_Model_tmax_effect},
we observe that the numerical solution gets close to, but fails to
converge to, the numerical steady state and, simultaneously, the time-step sizes stabilize to
a limiting step size $dt_\infty$.   (The figure only shows the behaviour up to
$t=1.5$, however the observed behaviour continues past this time for as long
as we chose to compute the solution.)  Figure \ref{Full_Model_tmax_effect}
is discussed fully in Section \ref{numerical_tests_PNP}.



When the VSSBDF2 adaptive time-stepper is taking (essentially) constant time steps it is effectively an SBDF2 time-stepper with
time-steps equal to  $dt_\infty$.  For for this reason, to try and understand the 
unexpected behaviour of the VSSBDF2 adaptive time-stepper,
we perform a stability 
analysis of the SBDF2 scheme linearized about the steady-state solution. 
 A significant challenge is that
the linearized scheme cannot be reduced via diagonalization to a study
of the scheme's behaviour for a single linear ODE.

We demonstrate
that the linearized scheme
is conditionally stable with a stability restriction $dt^*$.
Depending on the physical parameter values, when $dt = dt^*$, either there is one 
eigenvector with 
eigenvalue $-1$ or there is a pair of eigenvectors with complex eigenvalues of magnitude $1$.  The eigenvectors
are not highly oscillatory and, when performing the stability analysis using different spatial discretizations, we find
that the stability restriction, $dt^*$, does not change significantly.   Specifically, $dt^*$ does not go to zero
as $dx$ goes to zero; this shows that the conditional stability is not of
``diffusive type'' in which high frequencies can grow exponentially in time if $dt$ is too large.  

\begin{figure}[htb!]
\centering
\includegraphics[width=0.47\linewidth]{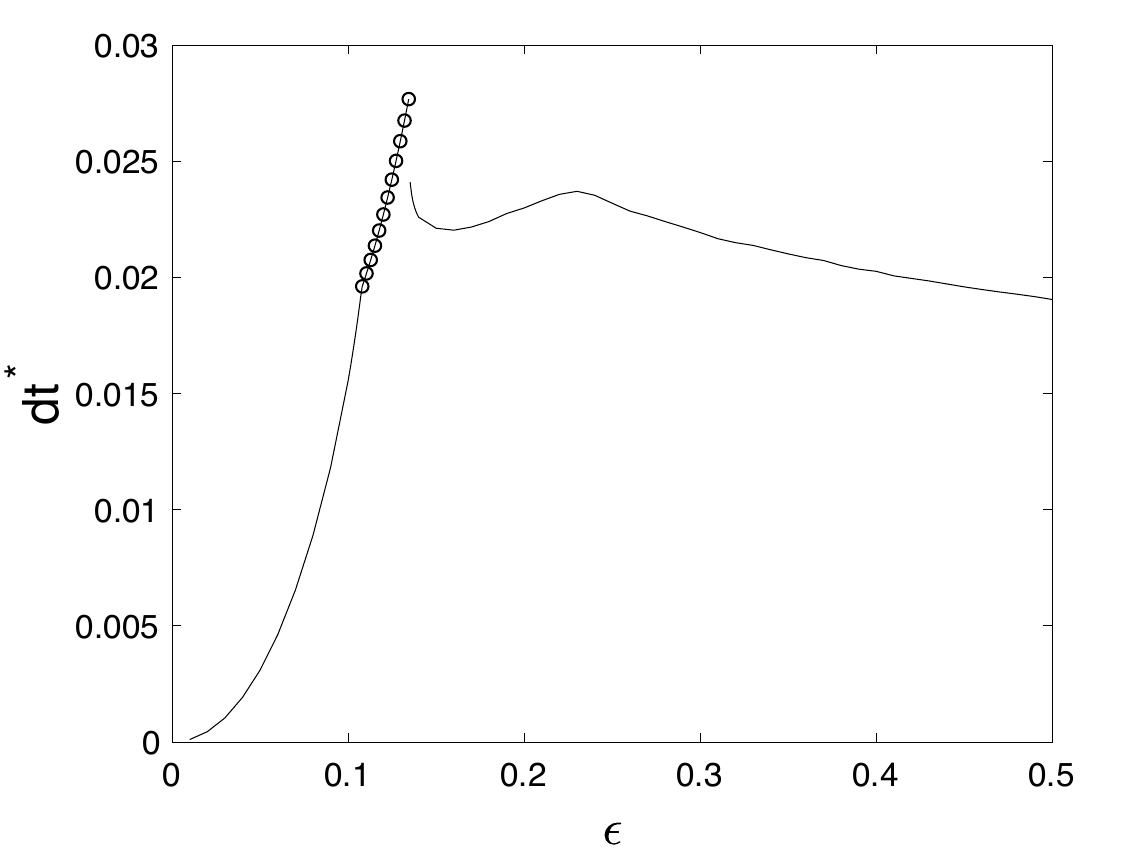}
\includegraphics[width=0.47\linewidth]{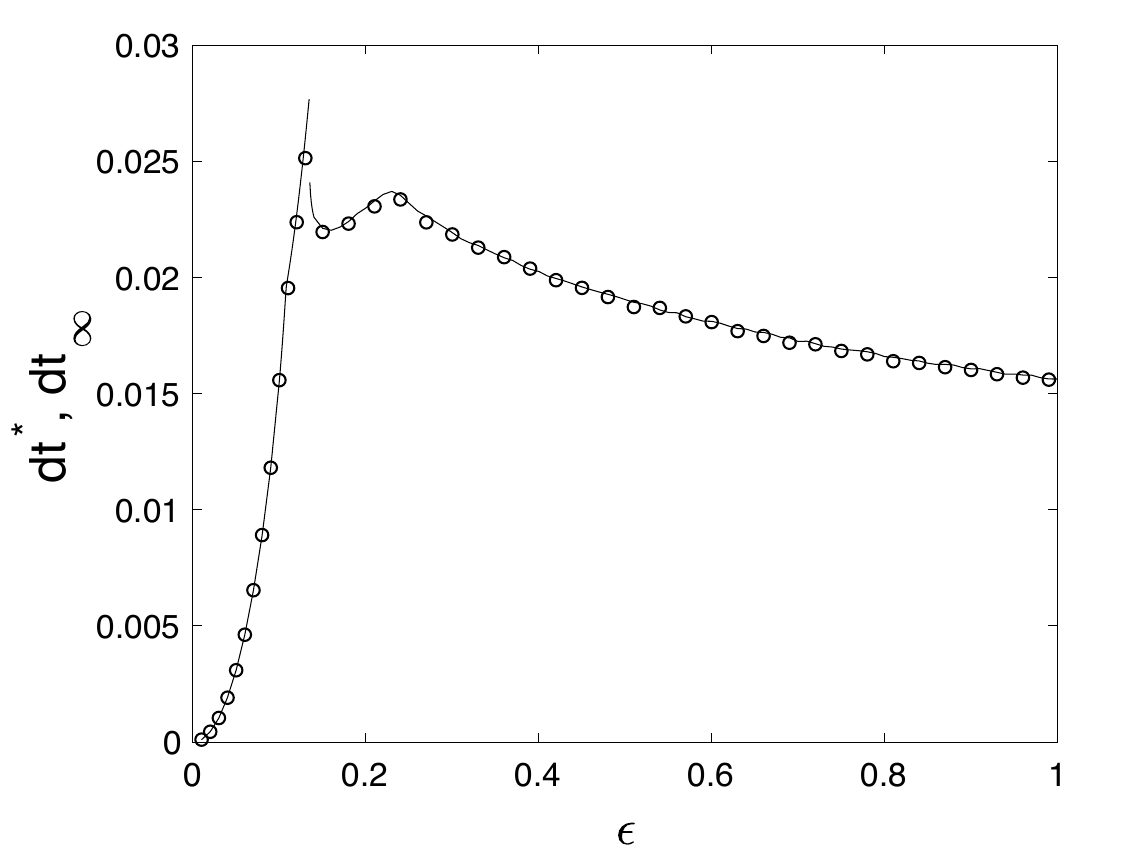}
\caption{PNP-FBV system \eqref{concentration_nondim}--\eqref{phi_bc_nondim_R} 
with physical and numerical parameters identical to those used for Figure \ref{Full_Model_tmax_effect} except
for $\epsilon$, which varies.
\underline{Left plot:} Solid line: $dt^*$ versus $\epsilon$ where $dt^*$ is found from the linear stability
analysis.  Open circles denote $\epsilon$ values for which the time-step stability restriction arises from a pair of complex conjugate eigenvalues
on the unit circle: $\epsilon \in (0.107764,0.134504)$.
\underline{Right plot:} 
Solid line: $dt^*$ versus $\epsilon$ where $dt^*$ is found from the linear stability
analysis.  Open circles: $dt_\infty$ as found by the time-step size stabilizing in the VSSBDF2 adaptive
time-stepper.  
}
\label{Full_model_dtthresh_vs_epsilon}
\end{figure}

The PNP-FBV system has a singular perturbation parameter $\epsilon$.
By varying $\epsilon$, we are able to study how the stability
restriction, $dt^*$, depends on $\epsilon$.  Figure \ref{Full_model_dtthresh_vs_epsilon} presents the graph of 
$(\epsilon,dt^*(\epsilon))$.  We
refer to the region below the
graph as the ``stability domain'' of the SBDF2 scheme linearized about the
state solution of the  PNP-FBV system.
We
find that the stability domain is not smooth --- there can be corners and jump
discontinuities in the graph of $dt^*$; see the left plot in Figure \ref{Full_model_dtthresh_vs_epsilon}.  Jump discontinuities are especially striking because they
mean that 
the same value of $dt$ could yield a stable SBDF2 computation for one value of $\epsilon$ but could result
in a computation that blows up for another, close value of $\epsilon$.   We have not seen this type of
phenomenon (non-smooth stability domains) reported in the literature.  We find that for small values of $\epsilon$ the stability domain
is not significantly influenced by the value of the (constant) imposed voltage or imposed current.

Our simulations suggest that $dt_\infty = dt^*$: the VSSBDF2 adaptive time-stepper is {\it finding} the
stability restriction for the SBDF2 scheme.   The right plot of Figure \ref{Full_model_dtthresh_vs_epsilon} presents the limiting time-step
found by the adaptive time-stepper ($dt_\infty$, open circles) as well as the
stability restriction ($dt^*$, solid line).  The Figure 
is fully discussed in Section \ref{linear_stability}.

The VSSBDF2 adaptive time-stepper can be used with, or without, a Richardson extrapolation step \eqref{final_extrapolation}.  Richardson
extrapolation is a common way to increase the accuracy of a scheme.  In
\cite{YPD_Part2_2019}, we demonstrate that Richardson extrapolation
can affect the linear stability of a scheme in various ways.  There, we give an example of a PDE for which the 
SBDF2 scheme is unconditionally
stable however, when used in combination with a Richardson extrapolation
step, the time-stepping scheme is conditionally stable.\\

\noindent {\bf Time-step Stability Restrictions}

\noindent
Time-step stability restrictions can arise in a variety of ways.
This article is not about stiffness in the sense of there being 
a restriction on the time-step size due to physical effects such as fast time scales.
Rather, it is about a system in which the fast time scales that arise from discretizing
the diffusion terms are well-handled by
a commonly-used semi-implicit scheme. However, as the solution equilibrates, the scheme becomes conditionally stable, leading to a time-step restriction.

There have been many approaches to the challenge of ``stiffness reduction'' whether in semi-implicit linear multi-step methods \cite{Rosales_Seibold_Shirokoff, Seibold_Shirokoff_Zhou, Akrivis_K_2003, Akrivis_KK_2003, Frank_1997, Douglas_Dupont_1971, Hou_Lowengrub_Shelley, Eggers_IMEX}, Runge-Kutta methods \cite{Abdulle_Medovikov_2001, Shin_Lee_Lee_2017, Izzo}, matrix exponential/integrating factor methods \cite{Kassam_Trefethen_2005, Ju_Zhang_Zhu_2015, Milewski_Tabak_1999}, and other approaches \cite{Minion_2003, Bruno_Jimenez_2014, Bruno_Lyon_2010a, Bruno_Lyon_2010b}.  These citations are provided as examples of the many seminal/well-written works in a large literature on the topic.  
It is
likely that some of the methods proposed in the cited works could allow one to use knowledge of the
structure of the linearized operator about the steady-state solution so as to modify the time-stepping scheme and remove the conditional stability.\\


%

\noindent {\bf Relevance of this work to other PDEs and other time-stepping schemes}

\noindent
In this article, we study the stability domain of the fixed-time-step SBDF2
scheme for the PNP-FBV system using a VSSBDF2 adaptive time-stepper.  We then confirm, and better understand, the findings by studying a
linearization of the SBDF2 scheme.  This two-pronged approach of using both
an adaptive time-stepper and a linearization study is not 
restricted to the PNP-FBV system
or to the VSSBDF2, SBDF2 pairing.
 
We used the VSSBDF2 adaptive time-stepper on several
dissipative systems that have non-constant, asymptotically stable steady states.
For example, we found that a basic reaction diffusion 
equation, $u_t = u_{xx} \pm u^2$,
did not yield a conditionally stable scheme when the SBDF2 scheme
is linearized around the asympotically stable steady state.
However, we did find some simple models related to the PNP-BDF system with asymptotically
stable steady states for which the SBDF2
scheme is conditionally stable when linearized about them \cite{YanThesis}.
The structure of the stability domain is also problem-specific.  We considered some other dissipative systems and 
did not find stability domains with corners, cusps, or jumps (otherwise we would have presented results for a simpler system
than the PNP-FBV system).

In terms of time-steppers,
if one is using a linear multistep method (LMM)
to study a physical system that has
asymptotically stable steady states, our approach is relevant in the
following ways.  
\begin{itemize}
\item By using an adaptive time-stepper that is built upon the variable step-size
version of the LMM, one is freed from needing to guess a ``good'' time-step size
for a fixed-time-step scheme.
Either the time-steps will get larger and larger as the solution converges to 
the steady state or they will stabilize to some limiting value $dt_\infty$.
In the latter case, if one
wishes to compute the numerical steady-state solution
up to round-off error, one then can do this with confidence by using the
constant-time-step scheme with a time-step size chosen smaller than
$dt_\infty$. 
\item We give a heuristic argument based on the local
truncation error as to why, in general, an adaptive time-stepper would
naturally find the stability restriction if the underlying
constant-time-step scheme is conditionally stable when linearized
about the steady state.  As a result, if one builds an adaptive
time-stepper based on a variable step-size version of the LMM
being used, one can use the adaptive time-stepper to explore the
stability domain of the constant-step-size LMM.  If, when computing an initial
value problem, the
time-step size stabilizes to a value $dt_\infty$, this suggests
that the underlying constant-time-step scheme is conditionally stable
when linearized about the steady state,
with stability restriction $dt^*$ and
$dt_\infty = dt^*$.
\item If the constant-step-size LMM is conditionally stable when
linearized about the steady state, one needs to
study the eigenvalues of the linearized
problem in order to to understand stability domain features such as
corners, jumps, and whether or not it is a single real-valued
eigenvalue that goes unstable as $dt$ exceeds $dt^*$.  The procedure
we use to linearize about the steady state and find the eigenvalues
and eigenvectors of the linearized system could
be used
for any LMM.  
\end{itemize}

%

\subsection{Structure of the article}
This article is structured as follows. Section \ref{PNP_system}
presents the PNP-FBV system.  Section \ref{crash_course} presents an
overview of the adaptive time-stepper.
Section \ref{numerical_tests_PNP} presents the
numerical simulations of the an initial value problem.
In Section
\ref{linear_stability}, the linearization about the steady-state
solution is presented.
In Subsection \ref{find_dt*}, the process for finding the stability 
restriction $dt^*$ is discussed.
The dependence of the stability
restriction $dt^*$ on physical quantities such as the singular
perturbation parameter $\epsilon$ and the imposed voltage are
discussed in Subsections \ref{eps_dependence} and \ref{v_dpendence}
respectively.
The dependence of the stability
restriction $dt^*$ on numerical aspects such as the mesh
and whether or not Richardson extrapolation is used
are
discussed in Subsections \ref{discretization_effect} and \ref{RE_effect}
respectively.

\section{The PNP-FBV system}
\label{PNP_system}

The Poisson-Nernst-Planck (PNP) equations describe the transport of charged species subject to diffusion and 
electromigration.   They have wide applicability in electrochemistry, and have been used to model a  
number of different systems, including porous 
media \cite{Biesheuvel2010PRE, Biesheuvel2011, Biesheuvel2012, Peters2016},
microelectrodes \cite{Streeter2008, Compton2011},
ion-exchange membranes \cite{Dydek2013, Nikonenko2010},
electrokinetic phenomena \cite{Yaroshchuk2012, Bazant2010, Bazant2009},
ionic liquids \cite{Bazant2011, Kornyshev2007},
electrochemical thin films \cite{Bazant2005, Chu2005, Biesheuvel2009},
fuel cells \cite{BiesheuvelFranco2009},
supercapacitors \cite{Lee2014},
and many more.
The Frumkin-Butler-Volmer boundary conditions describe charge transfer reactions at electrodes.

The one-dimensional, nondimensionalized  PNP equations 
for a medium with $2$ mobile species is
\begin{align}
\label{concentration_nondim}
\frac{\partial c_\pm}{\partial t} &= -\frac{\partial}{\partial x}\left[-\frac{\partial c_\pm}{\partial x} - z_\pm \, c_\pm \, \frac{\partial \phi}{\partial x}\right], \qquad t>0, \, x \in (0,1),\\
\label{poisson_nondim}
-\epsilon^2\frac{\partial^2 \phi}{\partial x^2} &= \frac{1}{2}\left(z_+ \, c_+ + z_- \, c_-\right), \qquad \qquad \qquad x \in (0,1),
\end{align}
where $c_\pm$ and $z_\pm$ are the concentration and charge number of the positive/negative ion, $\phi$ is the potential and $\epsilon$ is the ratio of the Debye screening length to the inter-electrode width $L$. 
This width is used
in the nondimensionalization of the original modelling equations \cite{Yan2017}; the
domain $(0,L)$ is rescaled to $(0,1)$.
We consider a model in which the anion and cation have a single charge ($z_\pm = \pm 1$) and
the anion has no charge-transfer reactions at the electrode: $c_-$ has no-flux boundary conditions:
\begin{equation}
\label{no_flux_BCs}
- \left( - \frac{\partial  c_-}{\partial x} +  \, c_- \, \frac{\partial \phi}{\partial x}\right)\bigg |_{x=0} = \left(- \frac{\partial  c_-}{\partial x} +  \, c_- \, \frac{\partial \phi}{\partial x} \right)\bigg |_{x=1}  = 0.
\end{equation}
The cation is assumed to have a reaction at the electrodes involving the transfer of one electron; this is modelled
using  
generalized Frumkin-Butler-Volmer (FBV) boundary conditions:
\begin{align}
\label{bv_nondim_L}
- \left( - \frac{\partial  c_+}{\partial x} -  \, c_+ \, \frac{\partial \phi}{\partial x}\right)\bigg |_{x=0} &= F(t) := 4k_{c,a} \,
c_+(0,t) \, e^{- 0.5 \; \Delta \phi_\text{left}} - 4 \, j_{ r,a} \, e^{0.5 \; \Delta \phi_\text{left}}, \\
\label{bv_nondim_R}
\left(- \frac{\partial  c_+}{\partial x} -  \, c_+ \, \frac{\partial \phi}{\partial x} \right)\bigg |_{x=1} &= G(t) := 4k_{c,c} \,
c_+(1,t) \, e^{-0.5 \; \Delta \phi_\text{right}} - 4 \, j_{r,c} \, e^{0.5 \; \Delta \phi_\text{right}},
\end{align}
where $k_{c,a}$, $k_{c,c}$, $j_{r,a}$, and $j_{r,c}$ are reaction rate parameters; the second part of the subscripts $a$ and $c$ refer to the anode and cathode, respectively . Equations \eqref{bv_nondim_L}--\eqref{bv_nondim_R} model the electrodeposition reaction
$ \ce{C^+ + e^- <=> M }$
where M represents the electrode material. 
The Stern layer is a compact layer of charge that occurs in the electrolyte next to an electrode surface
\cite{Bazant2013ACR, Soestbergen2012}; $\lambda_S$ denotes the effective
width of this layer.  
In equations \eqref{bv_nondim_L}--\eqref{bv_nondim_R}, $\Delta
\phi_\text{left}$ and $\Delta \phi_\text{right}$ refer to the
potential differences across the Stern layers that occur at the anode
and cathode respectively.  Specifically,
\begin{equation} \label{Delta_phi_is}
\Delta \phi_\text{left} = \phi_\text{anode} - \phi(0,t) = -\phi(0,t), 
\; \Delta \phi_\text{right} = \phi_\text{cathode} - \phi(1,t) = v(t)-\phi(1,t)
\end{equation}
where the potential at the anode has been set to zero and $v(t)$
denotes the potential at the cathode.  In addition, the Poisson
equation \eqref{poisson_nondim} uses a mixed (or Robin) boundary condition \cite{Bazant2005,
  Chu2005, Biesheuvel2009},
\begin{align}
\label{phi_bc_nondim_L}
- \epsilon \, \delta \; \frac{\partial \phi}{\partial x}\bigg |_{x=0} & = \Delta \phi_\text{left} := - \phi(0,t), \\
\label{phi_bc_nondim_R}
+ \epsilon \, \delta \; \frac{\partial \phi}{\partial x}\bigg |_{x=1} & =\Delta \phi_\text{right} := v(t) - \phi(1,t),
\end{align}
where $\delta=\lambda_S/L$.  Finally, there is an ODE which ensures conservation of electrical
current at the electrode \cite{Moya1995, Soestbergen2010},
\begin{equation}
\label{current_conservation_nondim}
-\frac{\epsilon^2}{2} \, \frac{d \; }{dt} \phi_x(1,t)
=j_\text{ext}(t) - \left[k_{c,c} \, c_+\left(1,t\right) \, e^{-0.5 \; \Delta \phi_\text{right}} - j_{r,c} \, e^{0.5 \; \Delta \phi_\text{right}}\right],
\end{equation}
where $j_\text{ext}(t)$ is the external current.  We refer
to the PNP equations with the generalized Frumkin-Butler-Volmer
boundary conditions as the PNP-FBV system.


The device is operated in two regimes --- either the current or the
voltage at the cathode is externally controlled\footnote{We use ``imposed'' as short-hand for externally controlled.}.  If the voltage at
the cathode, $v(t)$, is 
imposed then the the PNP-FBV
system \eqref{concentration_nondim}--\eqref{poisson_nondim} with
boundary conditions \eqref{no_flux_BCs}--\eqref{bv_nondim_R} and
\eqref{Delta_phi_is}--\eqref{phi_bc_nondim_R} are numerically solved,
determining $c_\pm$ and $\phi$.  The current is found {a
  postiori} using equation \eqref{current_conservation_nondim}.  If
the current, $j_\text{ext}(t)$, is 
imposed then the ODE
\eqref{current_conservation_nondim} is part of the PNP-FBV system and
is numerically solved along with the PDEs,
determining $c_\pm$, $\phi$, and $\phi_x(1,t)$ simultaneously.
The voltage $v(t)$ is then found {a postiori}. 

\section{The adaptive time-stepper}
\label{crash_course}

The companion article \cite{YPD_Part1} and \cite{YanThesis} present the numerical scheme in full:
spatial discretization, boundary conditions, splitting scheme, and error control.

The method of lines is used to discretize the parabolic PDEs \eqref{concentration_nondim}.
The spatial discretization reduces the parabolic-elliptic
system of PDEs to a differential-algebraic system of 	equations.
The system is handled using a splitting method: the ODEs are time-stepped,
 the system of algebraic equations is solved, the ODEs are
time-stepped again, and so forth.

The linear diffusion terms
in the parabolic PDEs \eqref{concentration_nondim} yield stiff linear terms in the ODEs.   
The terms in  \eqref{concentration_nondim} that model 
electromigration are nonlinear, yielding nonlinear terms in the ODEs.
The nonlinear terms make using implicit time-stepping methods unappealing.
Semi-implicit, or implicit-explicit schemes, are often used for stiff problems
as a way to avoid a fully-implicit treatment.

Consider the ODE
$u'= g(u) + f(u)$
where $g(u)$ is a stiff linear term and $f(u)$ is a nonlinear term.  
Given $u^{n-1}$ at time $t^{n-1}=t^n-dt$
and $u^n$ at time $t^n$, the SBDF2 scheme determines
$u^{n+1}$ at time $t^{n+1}=t^n+dt$ 
via
\begin{equation}
\text{SBDF2:} \hspace{.2in} \frac{1}{dt}\left(\frac{3}{2}u^{n+1} - 2 \, u^{n} + \frac{1}{2}u^{n-1}\right)
=g(u^{n+1}) + 2 \, f(u^n)- f(u^{n-1}), \label{sbdf2}
\end{equation}
where 
$u^n$ approximates
$u(t^n)$ (see, for example, \cite{Ascher1995}).  Our VSSBDF2 adaptive
time-stepper is based on a second-order variable step-size
semi-implicit backwards differencing formula, introduced by Wang
and Ruuth \cite{Wang2008}, as a generalization of the SBDF2 scheme:
\begin{align}\notag
&\text{VSSBDF2:} \hspace{.2in} \frac{1}{dt_\text{now}}\left(\frac{1+2\omega}{1+\omega}u^{n+1} - (1+\omega)u^{n} + \frac{\omega^2}{1+\omega}u^{n-1}\right)
\\
&\hspace{2in}=g(u^{n+1}) + (1+\omega)f(u^n)-\omega f(u^{n-1}), \label{vssbdf2}
\end{align}
where $\omega= dt_\text{now}/dt_\text{old}$, 
 $u^{n-1}$ is at time 
$t^{n-1}=t^n-dt_\text{old}$, and 
$u^{n+1}$ is at time 
$t^{n+1}=t^n+dt_\text{now}$.  Note that if $dt_\text{now} = dt_\text{old}$,
VSSBDF2 reduces to SBDF2.
%

The VSSBDF2 adaptive time-stepper is described in detail in the
companion article \cite{YPD_Part1}; see also \cite{YanThesis}.  The key idea is: if one
has computed the (approximate) solution up to the current time,
$(u^l,t^l)$ for $l = 0,\dots,n$, one can use these solutions and the
time-stepper to choose a new time $t^{n+1}$ so that the local
truncation error $| u^{n+1}-u(t^{n+1}) |$ is ``small but not too
small''.

One cannot know the local truncation error if one does not know the
(exact) solution $u(t^{n+1})$; in practice one needs an approximation
of the local trucation error.
We do this by performing a ``coarse" time step and a ``fine" time step to compute $u_c^{n+1}$ and $u_f^{n+1}$, respectively, and then using equation \eqref{final_error} to approximate the truncation error $\epsilon^{n+1}_c$.
\begin{equation}
\label{final_error}
\epsilon^{n+1}_c = \frac{8\left(dt_\text{old}+dt_\text{now}\right)}{7dt_\text{old}+5dt_\text{now}}\left(u^{n+1}_c-u^{n+1}_f\right)
 \approx u_c^{n+1}-u(t^{n+1}).
\end{equation}
If the error is acceptable\textcolor{blue}{,} we advance in time. If the error is
unacceptable\textcolor{blue}{,} we choose a new $dt_\text{now}$ and try again.  
If $dt_\text{now}$ has been accepted, we take $u^{n+1} = u^{n+1}_c$.
An algorithmic overview
is given in Appendix \ref{ATS}.
A detailed discussion of the adaptive time stepping and error control
schemes can be found in \cite{YanThesis,YPD_Part1}.   

This section describes the approach for an ODE; $u^n \in \R$.  It
generalizes immediately to a system of ODEs with $\u^n \in \R^N$.

\section{Simulations of the PNP-FBV system}
\label{numerical_tests_PNP}


Figure \ref{Full_Model_tmax_effect} presents a simulation of an  
initial value problem
for the PNP-FBV system
\eqref{concentration_nondim}--\eqref{phi_bc_nondim_R} with constant
imposed voltage.  The initial data is fixed, as are all the other
physical parameters.  Solutions are computed using the VSSBDF2
adaptive time-stepper.

The top plot in the figure demonstrates that, after a short
transient, the solution initially
decays exponentially quickly to a numerical steady state.  However, once the
solution is within (approximately) $10^{-7}$ of the steady-state
solution, this convergence ends and the computed solution stays about
$10^{-7}$ away from the steady state.  The middle plot in the figure 
demonstrates
that the VSSBDF2 adaptive time-stepper is keeping the (approximate) local
truncation error \eqref{final_error}
within the user-specified interval.  
The bottom plot in the figure 
demonstrates that the time-step size initially increases exponentially
fast and after a while
it decreases and stabilizes to $dt_\infty$.  The dashed line in the bottom figure is
the stability restriction found by the linear stability analysis
discussed in Section \ref{linear_stability}:
$dt^*$.  This simulation demonstrates that the VSSBDF2 adaptive time-stepper
appears to eventually stabilize at a time-step size that is precisely
the stability restriction.

The top plot presents the logarithm of the norm of the deviation of  $\c_{+}$ from the numerical steady state $\c_{+,ss}$; that is $\log(\| \c_+^n - \c_{+,ss}\|)$.  The deviations of
$\c_-^n$ and $\pmb{\phi}^n$ from the respective numerical steady states behave similarly.
The numerical steady state, $\c_{\pm,ss}$ and
$\pmb{\phi}_{ss}$, satisfies the discretized version of the
steady-state equations $0 = c_{\pm,xx} +
z_\pm \, ( c_\pm \, \phi_x )$ and \eqref{poisson_nondim}.  To find them,
a simulation using the VSSBDF2 adaptive time-stepper is stopped
once the time steps have stabilized.  The SBDF2 time-stepper is then
used to continue the simulation with a (fixed) time-step $dt$ chosen to be
smaller than $dt_\infty$.  The local truncation error tends to zero exponentially
fast and the simulation is stopped once the computed solution
satisfies the discretized steady-state equations (up to round-off).  
This late-time
solution is taken as the numerical steady-state solution.

If, rather than studying
the
deviations, one approximates the time derivative, $\c_{+,t}$, using \eqref{vssbdf2},
then the plot of $\log(\|\c_{+,t}\|)$ versus $t$ will show that 
$\| \c_{+,t}\|$ decreases and then stabilizes at a {\it nonzero} value.

Two of the user-specified parameters
of an adaptive time-stepper are $dt_\text{min}$ and $dt_\text{max}$.
The time-stepper is not allowed to take $dt$ smaller than
$dt_\text{min}$ or larger than $dt_\text{max}$.  The above-described behaviour is
what is observed if $dt_\infty < dt_\text{max}$.  If, by
chance, it happens that $dt_\text{max}$ is smaller than $dt_\infty$ then as
the solution equilibrates the time steps increase to $dt_\text{max}$
and are then held at that value.  The solution subsequently converges to 
the steady-state solution.

Rosam, Jimack and Mullis \cite{Rosam2007} used an
adaptive SBDF2 algorithm to study a problem in binary alloy solidification. In their Figure 4, they appear to show time-steps
stabilizing to a constant value, but the reason is
not given: they report that it is related to the tolerance set in the adaptive time-stepper.
We did not observe such a phenomenon when we varied $tol$; we found the same
limiting time-step size $dt_\infty$.  Also, we find that decreasing $tol$ leads to the
deviations becoming smaller before levelling out (see top plot of Figure \ref{Full_Model_tmax_effect}).

\section{Numerical Linear Stability}
\label{linear_stability}

If the imposed voltage or imposed current is constant for a period of
time and if $dt_\text{max}$ is large, then the solution attempts to equilibrate and the VSSBDF2 adaptive
time-stepper stabilizes to take (nearly) constant time
steps $dt_\infty$.  The VSSBDF2 scheme with constant time steps is the
SBDF2 scheme; for this reason we study the SBDF2 time-stepper to try
and understand this stabilization of the VSSBDF2 adaptive
time-stepper.

Consider the
SBDF2 scheme \eqref{sbdf2} applied to the $N$ ODEs 
$\u_t = \f(\u)+\g(\u)$
that arise from spatially discretizing the PDE $u_t = f(u,u_x,u_{xx},\dots) + g(u,u_x,u_{xx},\dots)$: 
\begin{equation} \label{BDF2_IMEX}
\frac{1}{dt} \left( \frac{3}{2} \u^{n+1} - 2 \u^n + \frac{1}{2} \u^{n-1} \right)
= \g(\u^{n+1}) + 2 \f(\u^n) - \f(\u^{n-1}).
\end{equation}
The bold-faced quantities are vectors in $\R^N$.
A steady state satisfies $0 = \g(\u_{ss}) + \f(\u_{ss})$.   
Linearizing about $\u_{ss}$ yields 
\begin{equation} \label{vectorized_2}
\frac{1}{dt} \left( \frac{3}{2} \d^{n+1} - 2 \d^n + \frac{1}{2} \d^{n-1} \right)
= \mathbf{J}_{\g}(\u_{ss}) \, \d^{n+1} + 2 \, \mathbf{J}_{\f}(\u_{ss})  \, \d^n - \mathbf{J}_{\f}(\u_{ss})  \, \d^{n-1}.
\end{equation} 
where $\d^m = \u^m - \u_{ss}$ is the deviation from the steady state and $\mathbf{J}_{\f}(\u_{ss})$ and $\mathbf{J}_{\g}(\u_{ss})$ are the Jacobian matrices evaluated at $\u_{ss}$; e.g. $(\mathbf{J}_{\f}(\u_{ss}))_{ij} = \frac{ \partial f_i}{\partial u_j}(\u_{ss})$.   
For simple problems, $\mathbf{J}_{\f}$ and $\mathbf{J}_{\g}$ can be determined analytically and evaluated at $\u_{ss}$.
Otherwise, one can numerically approximate $\mathbf{J}_{\f}(\u_{ss})$ and $\mathbf{J}_{\g}(\u_{ss})$
in a variety of ways.  We used a simple centre-difference scheme.  

%

If the Jacobian matrices, $\mathbf{J}_{\f}(\u_{ss})$ and
$\mathbf{J}_{\g}(\u_{ss})$, can be simultaneously diagonalized, then
the system \eqref{vectorized_2} reduces to a decoupled system of
second-order difference equations.  In this case, the stability
analysis is straight-forward: one computes the roots of the
now-decoupled difference equations and analytically studies how they depend on $dt$ and the
eigenvalues of $\mathbf{J}_{\f}(\u_{ss})$ and
$\mathbf{J}_{\g}(\u_{ss})$.  Indeed, in Appendix A of
\cite{YPD_Part2_2019}, we present the linear stability
analysis of \eqref{BDF2_IMEX} for a single linear ODE.  We use it 
to study the logistic equation, demonstrating that that the $dt_\infty$ found by the VSSBDF2 adaptive time-stepper  is in sharp agreement with the
analytically-determined stability restriction $dt^*$.

For the PNP-FBV system \eqref{concentration_nondim}--\eqref{phi_bc_nondim_R}, 
the Jacobians {\it cannot} be simultaneously diagonalized.  
For this reason,
we proceed with a numerical computation of the eigenvalues and eigenvectors of the linearized scheme \eqref{vectorized_2} rewritten as
\begin{equation} \label{vectorized_3}
\mathbf d^{n+1} = M_\text{new}  \, M_\text{now}  \, \mathbf d^n +  M_\text{new}  \, M_\text{old}  \, \mathbf d^{n-1},
\end{equation}
where
$$
M_\text{new} = \left(\frac{3}{2} \, I - dt  \, \mathbf{J}_{\g}(\u_{ss}) \right)^{-1}, 
M_\text{now} = 2  \, I + 2 \, dt  \, \mathbf{J}_{\f}(\u_{ss}),  
M_\text{old} = -\frac{1}{2}  \, I - dt  \, \mathbf{J}_{\f}(\u_{ss}).
$$
Equation (\ref{vectorized_3}) is a system of $N$ second-order linear difference equations.  Solving it requires the initial deviation,
$\mathbf d^0$, as well as the deviation after one time-step, $\mathbf d^1$.  The system is rewritten 
\cite[\S D.2.1]{LeVeque_FD_book}
as $2N$ first-order
linear difference equations in the standard manner: $\mathbf D^n := [ \mathbf d^{n-1} \, ; \, \mathbf d^n] \in \R^{2N}$ and
$A$ is the companion matrix for the difference equation:
\begin{equation} \label{vectorized_4}
\mathbf D^{n+1} = A \, \mathbf D^n = 
\left( \begin{array}{cc}
0 & I \\
M_\text{new}  \, M_\text{old}   & M_\text{new}  \, M_\text{now} \end{array} \right) \, \mathbf D^n.
\end{equation}
If $(\lambda_j,\v_j)$ is an eigenvalue-eigenvector pair of $A$ then the structure of $A$ implies that
$\v_j =  [\mathbf d_j \, ; \, \lambda_j  \mathbf d_j  ]$ for some $\mathbf{d}_j \in \R^N$.  
If $A$ has $2N$ linearly independent eigenvectors, it follows that
the general solution of the linearized problem (\ref{vectorized_2}) is
\begin{equation} \label{general_solution}
\mathbf d^n = \sum_{j=1}^{2N} c_j \, \lambda_j^n \, \mathbf d_j
\end{equation}
where the $2N$ coefficients, $c_j$, are determined using $\mathbf d^0, \mathbf{d}^1 \in \R^N$.

The connection between the linearized scheme \eqref{vectorized_2} and the nonlinear scheme \eqref{BDF2_IMEX} 
is via the stability theory of fixed points for discrete dynamical systems.  If $\mathbf{J}_{\f}$ and 
$\mathbf{J}_{\g}$ are continuous in a neighbourhood of $\u_{ss}$ and if $dt$ is 
such that $(\frac{3}{2} \, I - dt  \, \mathbf{J}_{\g}(\u_{ss}))$ is invertible then the discrete dynamical system
\begin{equation} \label{DDS}
\U^{n+1} 
= \F(\U^n) :=
\begin{pmatrix}
\U^n_2\\
\left(\frac{3}{2} \, I - dt \, \g \right)^{-1}(2 \U^n_2 + 2 \, dt \, \f(\U^n_2) - \frac{1}{2} \U^n_1 - dt \, \f(\U^n_1))
\end{pmatrix}
\end{equation}
is defined in a neighbourhood of the fixed point $[\u_{ss} \, ; \, \u_{ss}]$.  Defining $\U^{n+1} = [ \u^n \, ; \, \u^{n+1}]$,
this discrete dynamical system
\eqref{DDS} is equivalent to the SBDF2 time-stepping scheme \eqref{BDF2_IMEX}.  The companion matrix
$A$ is the linearization of \eqref{DDS} at the fixed point  $[\u_{ss} \, ; \, \u_{ss}]$.  Therefore,
if all eigenvalues of $A$ have 
magnitude less than $1$, then $[\u_{ss} \, ; \, \u_{ss}]$ is an asymptotically stable fixed point of \eqref{DDS}
and $\u_{ss}$ is an asymptotically stable fixed point of the SBDF2 scheme \eqref{BDF2_IMEX}.

\subsection{Finding the stability restriction $dt^*$ (if there is one)}
\label{find_dt*}


Given a
particular PDE (or set of PDEs) and boundary conditions, we compute
the Jacobians (using the time-step $dt$ that was used in the SBDF2 time-stepper
to find $u_{ss}$) about the steady state and construct the matrix $A$ in
\eqref{vectorized_4}.
%
The eigenvalues of $A$ will have magnitude
less than one --- otherwise the time-stepper would not have found the steady-state solution.
%
%
To determine if there is a linear stability restriction, one
increases the
time-step size, recomputes A and its eigenvalues, and then determines
if any eigenvalues have magnitude greater than one for this new value
of $dt$.
Proceeding in this way, one seeks a time-step size at which an
eigenvalue(s) crosses from magnitude less than one to magnitude greater
than one.  We use an iterative bisection method to approximate this critical time-step size, $dt^*$.   The scheme is conditionally stable in that $\u_{ss}$ is
an asymptotically stable steady-state solution of the SBDF2 time-stepper if $dt < dt^*$ and is an unstable steady state
if $dt > dt^*$.  

Whether or not there is such a stability restriction on $dt$ is a problem-specific
question.  And, of course, not finding $dt^*$ could either indicate unconditional stability or it 
could indicate that one has not tried large enough values of $dt$.

The above process depends on first finding a value of $dt$ 
such that $u_{ss}$ is
asymptotically stable under the dynamics \eqref{BDF2_IMEX}.
In practice, $dt^*$ is not a priori known, and so the region $0 < dt < dt^*$
is unknown.
For this reason, we first use the VSSBDF2 adaptive 
time-stepper. If it yields a solution that appears to be ``trying to but failing to
to converge''\footnote{That is, the local truncation error
stays in the user-specified interval, the approximation of $u_t$ decays exponentially and then stabilizes, and 
the time-step sizes are stabilizing at some value $dt_\infty < dt_\text{max}$}
to a steady-state solution,
we take this as evidence of the SBDF2 scheme's being conditionally stable
when linearized about the steady state.  
We then choose some $dt < dt_\infty$ and repeat the simulation
using the (constant time-step) SBDF2 scheme to find a steady-state solution,
as described in Section \ref{numerical_tests_PNP}.  \\

Turning to the PNP-FBV system \eqref{concentration_nondim}--\eqref{phi_bc_nondim_R},
for a fixed constant imposed voltage and $\epsilon$, we use the (constant time-step) SBDF2 
time-stepper to find the discrete
steady-state solution $\c_{+,ss}$, $\c_{-,ss}$, and $\pmb{\phi}_{ss}$.  
The steady-state concentrations
are
concatenated into one vector $\u_{ss} := [\c_{+,ss} \, ; \, \c_{-,ss}]$.  The right-hand sides of  the discretized
evolution equations \eqref{concentration_nondim}
are similarly concatenated: $\f$ is the spatial discretization of $[(c_+ \, \phi_x)_x \, ; \, (c_- \, \phi_x)_x]$
and $\g$ is the discretization of $[c_{+,xx} \, ; \, c_{-,xx}]$.   
We then approximate the Jacobians $\mathbf{J}_{\f}(\u_{ss})$ and $\mathbf{J}_{\g}(\u_{ss})$.
If there are $N$ mesh points then $\u_{ss} \in \R^{2N}$ and 
the Jacobians are $2N \times 2N$ matrices.   A value of $dt$ is chosen and the 
$4N \times 4N$ companion matrix $A$ in (\ref{vectorized_4}) is constructed and its eigenvalues and eigenvectors
are computed.   The value of $dt$ is then increased and the process is repeated.

\begin{figure}[htb!]
\centering
\includegraphics[width=0.45\linewidth]{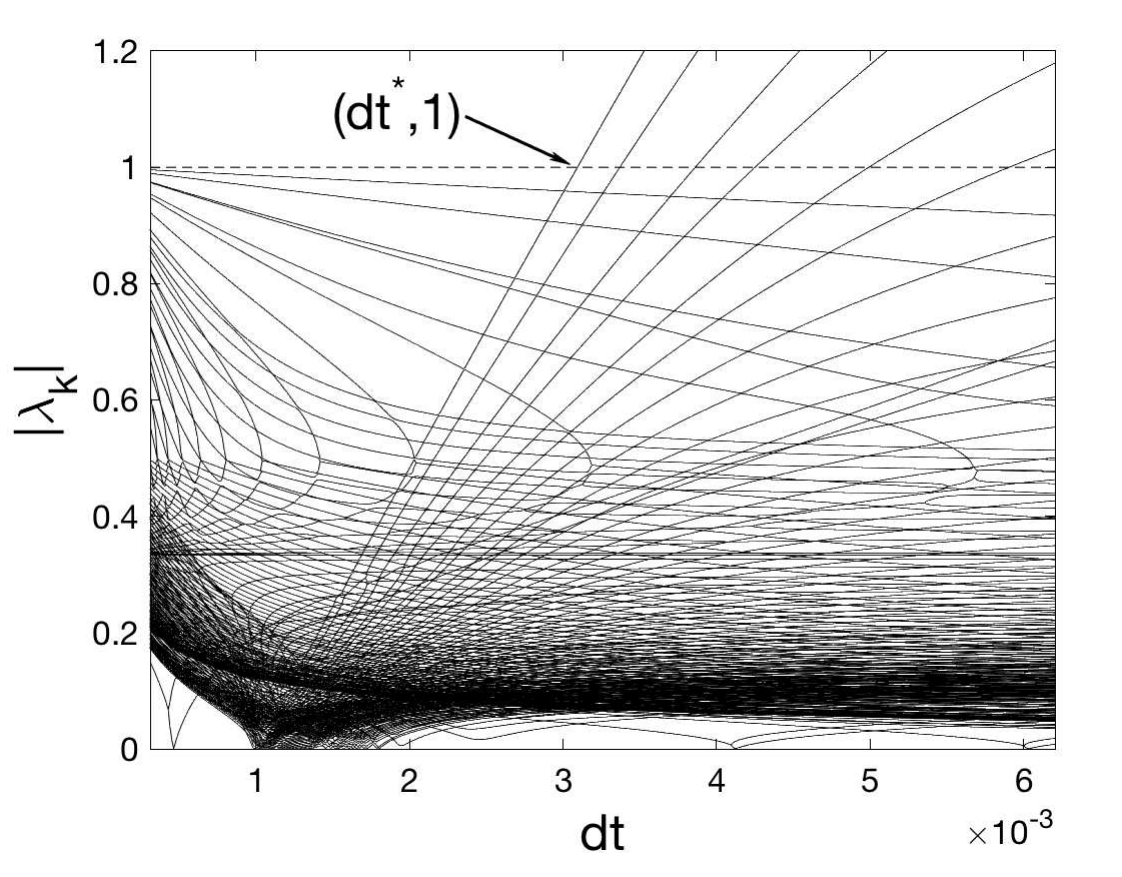}
\includegraphics[width=0.45\linewidth]{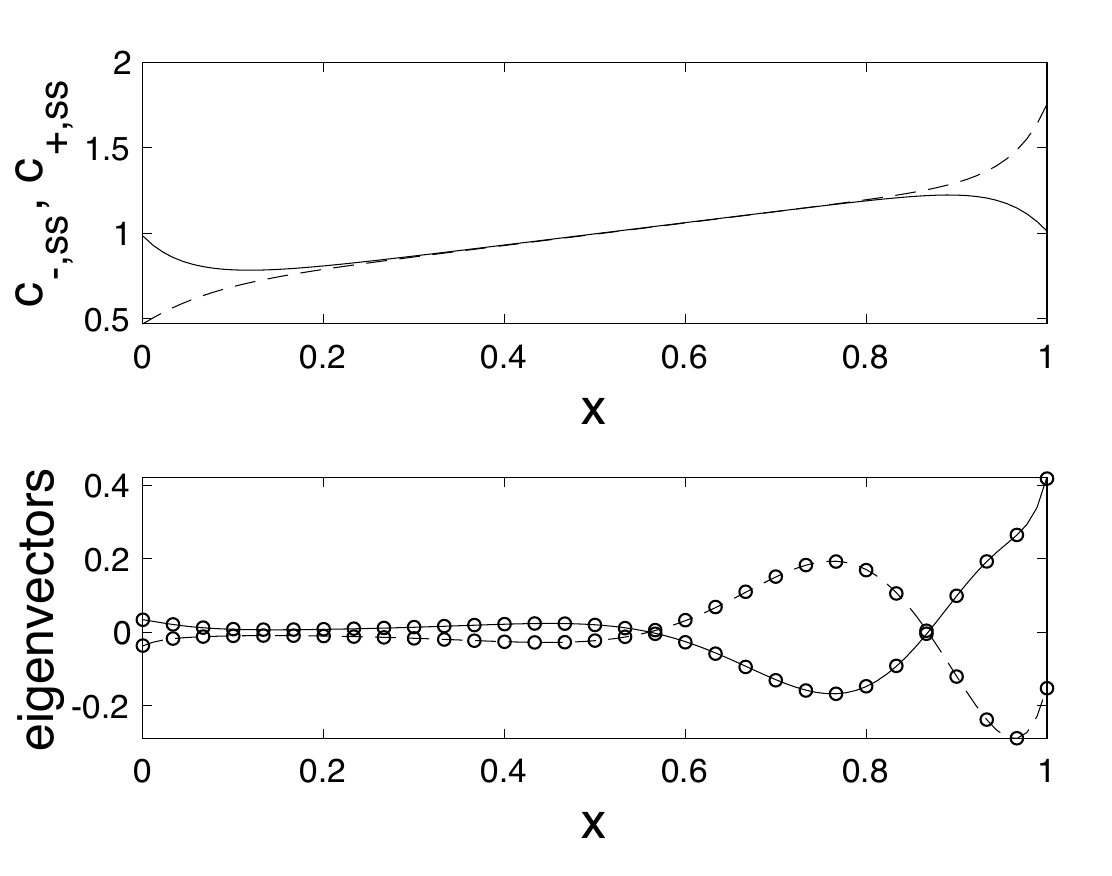}
\caption{PNP-FBV system \eqref{concentration_nondim}--\eqref{phi_bc_nondim_R}
with physical and numerical parameters identical to those used for Figure \ref{Full_Model_tmax_effect}. 
\underline{Left plot:} The magnitudes of the 364
eigenvalues are plotted versus $dt$.  
The largest magnitude branch crosses at $dt^*= .003094$.
\underline{Right plot, top:} The steady-state profiles --- the solid line is
$\c_{+,ss}$ and the dashed line is $\c_{-,ss}$.  
\underline{Right plot, bottom:} At $dt = dt^*$, one real-valued eigenvalue is ``barely stable'': $\lambda = -1$. The corresponding ``barely stable'' eigenvectors are plotted --- the solid line is
the eigenvector for $\c_+$ and the dashed line is the eigenvector for $\c_-$. 
Both eigenvectors have been chosen to have $l^2$ norm $1$.  The open circles denote late-time deviations
from the steady states, as computed using the VSSBDF2 adaptive
time-stepper.  The deviations have been normalized to have $l^2$ norm $1$; only a third of the  $N=91$ data
points are plotted for tidiness.  
}
\label{Full_Model_eps_0p05}
\end{figure}

Figure \ref{Full_Model_eps_0p05} presents results for $\epsilon = .05$ with 
constant imposed voltage.
In the left figure, the magnitudes of all eigenvalues are plotted --- we see that for small
values of $dt$, all eigenvalues have magnitude less than one and that, as $dt$ is increased, one branch
goes unstable.  We follow this branch to find the time-step size at which the magnitude 
equals 1; this is the stability restriction $dt^*$.  
We find that one eigenvalue
crosses the unit circle, crossing at value $-1$.  In the top plot in the right figure, we plot the steady
states $\c_{\pm,ss}$.  In the bottom plot in the right figure, we plot the eigenvectors at the stability
restriction $dt^*$.  We refer to the eigenvector-eigenvalue pair as ``barely stable''.

To demonstrate that the ``barely stable''
eigenvalue-eigenvector pair is the cause of the failure to converge to
the steady state shown in the top plot of Figure \ref{Full_Model_tmax_effect}, we took $\c_+$ and $\c_-$ at a late time ($t=100$) and computed the
corresponding deviations from the steady state $\d_+$ and $\d_-$.  
In
the bottom-right plot of Figure \ref{Full_Model_eps_0p05}, the normalized deviations are plotted with open
circles --- note that they closely fit the ``barely stable'' eigenfunctions.\\

To see why it is unsurprising that the VSSBDF2 adaptive time-stepper would 
adjust its timesteps until they stabilize
at the stability restriction of the underlying SBDF2 scheme, we consider the local truncation error for
the SBDF2 scheme, applied to the ODE $u_t = f(u)+g(u)$, close to a steady state $u_{ss}$:
\begin{align*}
LTE &= u^{n+1}-u(t_{n+1}) = d^{n+1}-d(t_{n+1}) \\
& =\left(\frac{2}{3} d^{\prime \prime \prime}(t_n) - g''(u(t_n)) \, d'(t_n)^2- g'(u(t_n)) \, d''(t_n) \right)dt^3 + O(dt^4).
\end{align*}
Here $u(t_n) = u^n$, $d(t_{n+1}) = u(t_{n+1})-u_{ss}$, and $d^{n+1}=u^{n+1}-u_{ss}$.  
Assume the dynamics of the underlying system of ODES is that of solutions converging
to an asymptotically stable steady state.  Now assume that that the numerical
approximations are also converging to the asymptotically stable steady state.
In this situation,
the deviation is decaying exponentially in time: $d(t) \cong C \exp(-\lambda t)$
and the local truncation error can be bounded
\begin{equation} \label{LTE_bound}
|C| \, \alpha \, e^{-\lambda t_n} dt^3 
\leq  |LTE| \leq |C| \, \beta \, e^{-\lambda t_n} dt^3 
\end{equation}
where $\alpha$ and $\beta$ are determined by $\lambda$ and uniform bounds on $g'$ and $g''$ near $u_{ss}$.
The upper bound in \eqref{LTE_bound} implies that if $dt$ is held fixed, the LTE will decay to zero as $t_n \to \infty$.  The lower bound in \eqref{LTE_bound} implies that if the LTE is required to satisfy
a constraint such as $LTE \geq tol-range > 0$, then $dt$ must grow exponentially as $t_n \to \infty$.

Whether or not the the numerical
approximations are converging to the asymptotically stable steady state
is determined by the spectral radius of the scheme linearized about
the steady state.
For a system of $N$ ODEs, the spectral radius
of the linearized scheme is
$$
|\lambda(dt)|_\text{max} = \max_{1 \leq i \leq N}\{ |\lambda_i(dt)| \}.
$$
If $dt$ is such that $|\lambda(dt)|_{\text{max}}<1$,
then the LTE for the SBDF2 scheme
will go to zero exponentially fast as the number of time steps goes
to infinity: the numerical
approximations are converging to the asymptotically stable steady state.
Similarly, if $|\lambda(dt)|_{\text{max}}>1$, then the LTE will grow exponentially
until nonlinear effects become relevant.
When our simulations with the VSSBDF2 adaptive time-stepper stabilize
to $dt_\infty$, they 
are taking essentially-constant time-steps.  At the same time, the LTE is constrained to stay in an interval $[tol-range,tol+range]$ where $tol-range>0$.  This can only happen if $|\lambda(dt_\infty)|_{\text{max}}=1$; i.e., $dt_\infty = dt^*$.

There is nothing specific to the SBDF2 scheme and its variable step-size generalization VSSBDF2 in this argument.  Any adaptive time-stepper that
is built on a variable step-size generalization of a constant-step-size scheme
could be used to explore the stability properties of the constant-step-size scheme.

\subsection{Dependence of the stability domain on the singular perturbation parameter $\epsilon$}
\label{eps_dependence}

\begin{figure}[htb!]
\centering
\includegraphics[width=0.47\linewidth]{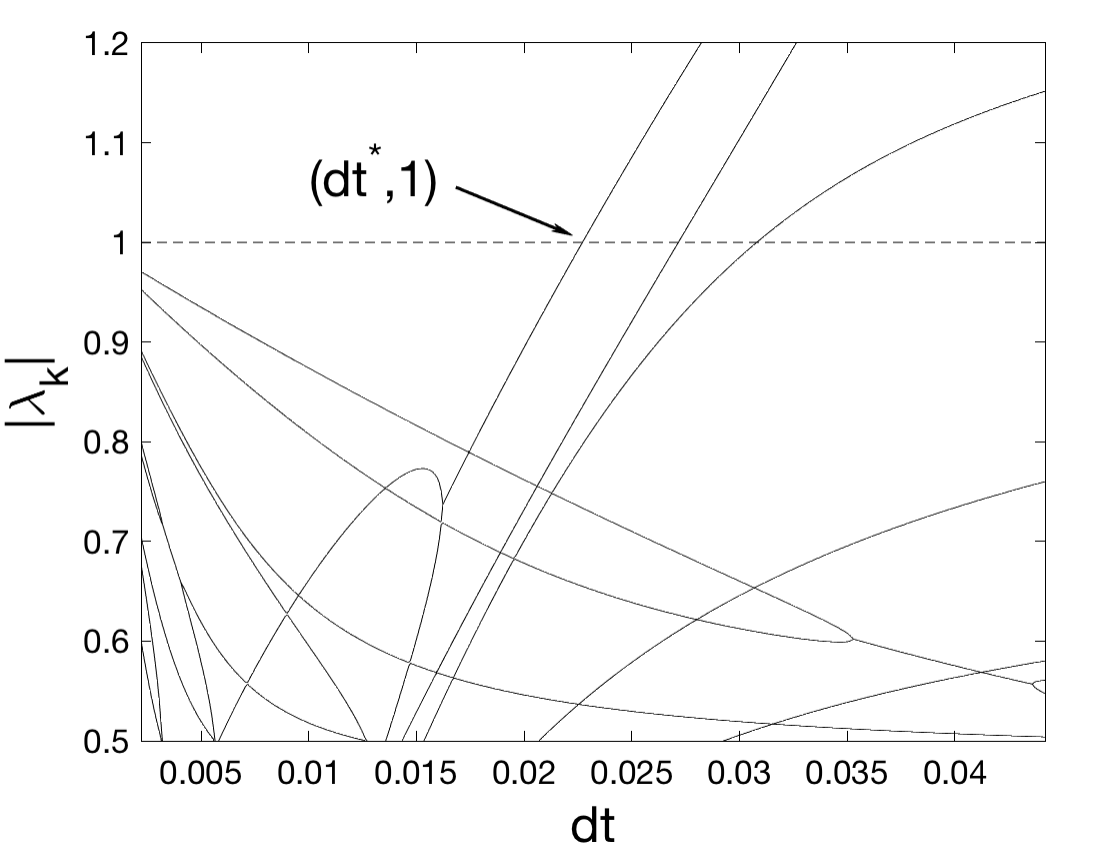}
\includegraphics[width=0.47\linewidth]{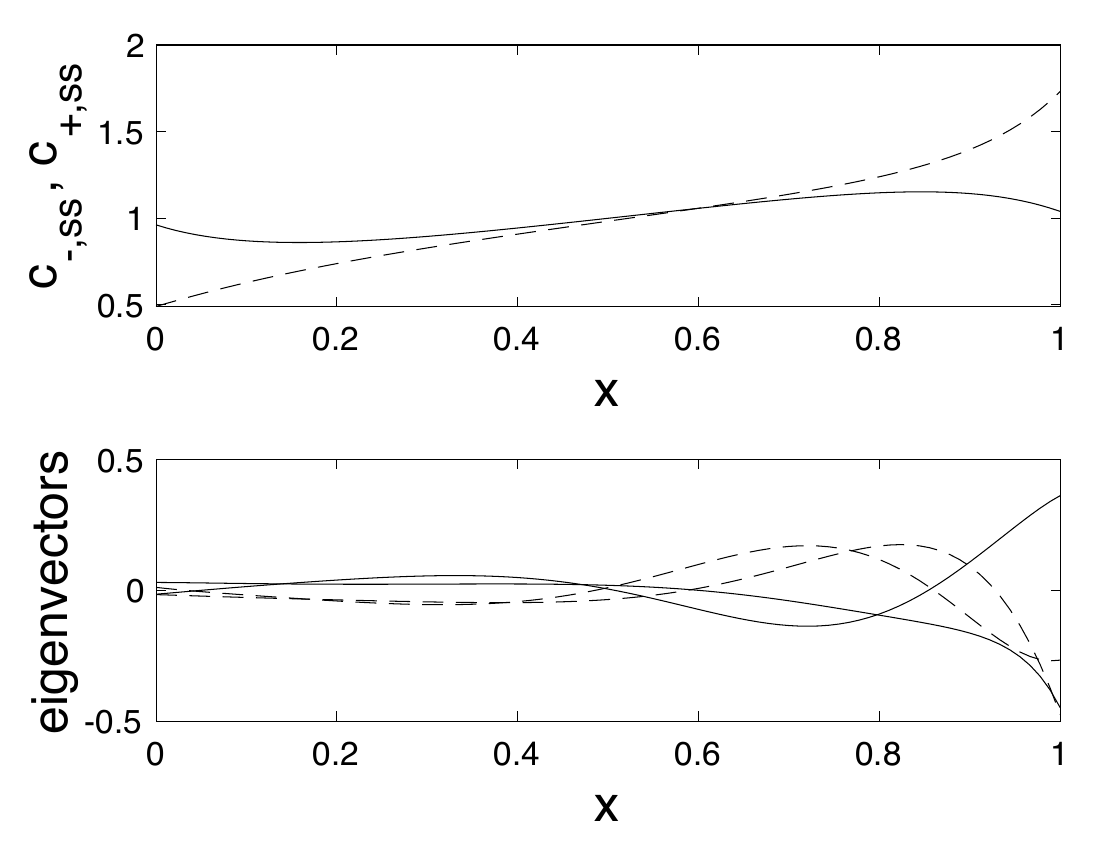}
\caption{PNP-FBV system \eqref{concentration_nondim}--\eqref{phi_bc_nondim_R} 
with physical and numerical parameters identical to those used for Figure \ref{Full_Model_tmax_effect} except
for $\epsilon = .12$.
\underline{Left plot:} The magnitudes of the 364
eigenvalues are plotted versus $dt$; the vertical range has been truncated for a tidier plot.  
The largest magnitude branch crosses at $dt^*= .02271$.
\underline{Right plot, top:} The steady-state profiles --- solid line is
$\c_{+,ss}$ and dashed line is $\c_{-,ss}$.  
\underline{Right plot, bottom:} A pair
of complex-valued eigenvalues, $-0.9797 \pm 0.2008 \, i$, go 
unstable.
The corresponding eigenvectors are plotted --- the solid lines are
the unstable eigenvectors for $\c_+$ and the dashed line are the unstable eigenvectors for $\c_-$. }\label{Full_Model_eps_0p12}
\end{figure}

We now consider the stability properties of the PNP-FBV system \eqref{concentration_nondim}--\eqref{phi_bc_nondim_R}
for a range of values of $\epsilon$, holding the imposed voltage and all other parameters fixed..  We find that for $\epsilon \in (0.107764,0.134504)$ the instability takes the form of a pair of complex
eigenvalue crossing the unit circle; for all other values we considered it was a single eigenvalue crossing
at $-1$.  Figure \ref{Full_Model_eps_0p12} is the analogue of Figure \ref{Full_Model_eps_0p05} but
for a value of $\epsilon$ that results in two complex eigenvalues crossing the unit circle.  

For each $\epsilon$, 
we find the stability restriction $dt^*$.  The left plot of Figure
\ref{Full_model_dtthresh_vs_epsilon} presents $dt^*$ as a
function of $\epsilon$.   The open circles indicate the interval of $\epsilon$
values for which a pair of complex conjugate eigenvalues cross the unit circle.
The stability restriction $dt^*$
is a continuous function of $\epsilon$ except for a jump at $\epsilon \approx 0.134504$.
Also, $dt^*$ appears to be a smooth function of $\epsilon$ except at
$\epsilon \approx 0.134504$ (where there is a jump in $dt^*$) and at 
$\epsilon \approx 0.107764$ (where there is a corner).

The jump in $dt^*$ is striking --- if one were using an SBDF2 time-stepper with 
$dt = .025$ then this would yield a stable
simulation for $\epsilon$ which is close to, but slightly smaller than, the critical value 
of $\epsilon \approx 0.134504$.  The simulation would
be unstable simulation for $\epsilon$ which is close to, but slightly larger than,
this critical value.  The stability of the SBDF2 simulation is not a continuous function of the 
parameter $\epsilon$.

The right plot of Figure
\ref{Full_model_dtthresh_vs_epsilon} compares $dt_\infty$ as found from
the VSSBDF2 adaptive time-stepper to $dt^*$ as found from the linear
stability study of the steady state.  The solid line plots
$dt^*$ versus $\epsilon$; the circles plot $dt_\infty$.  The
circles align closely with the solid lines, providing
compelling evidence that it is the numerical instability of the scheme
near the steady state which is causing the VSSBDF2 adaptive time-stepper to
stabilize its time-steps.  

\begin{figure}[htb!]
\centering
\includegraphics[width=0.47\linewidth,height=0.448\linewidth]{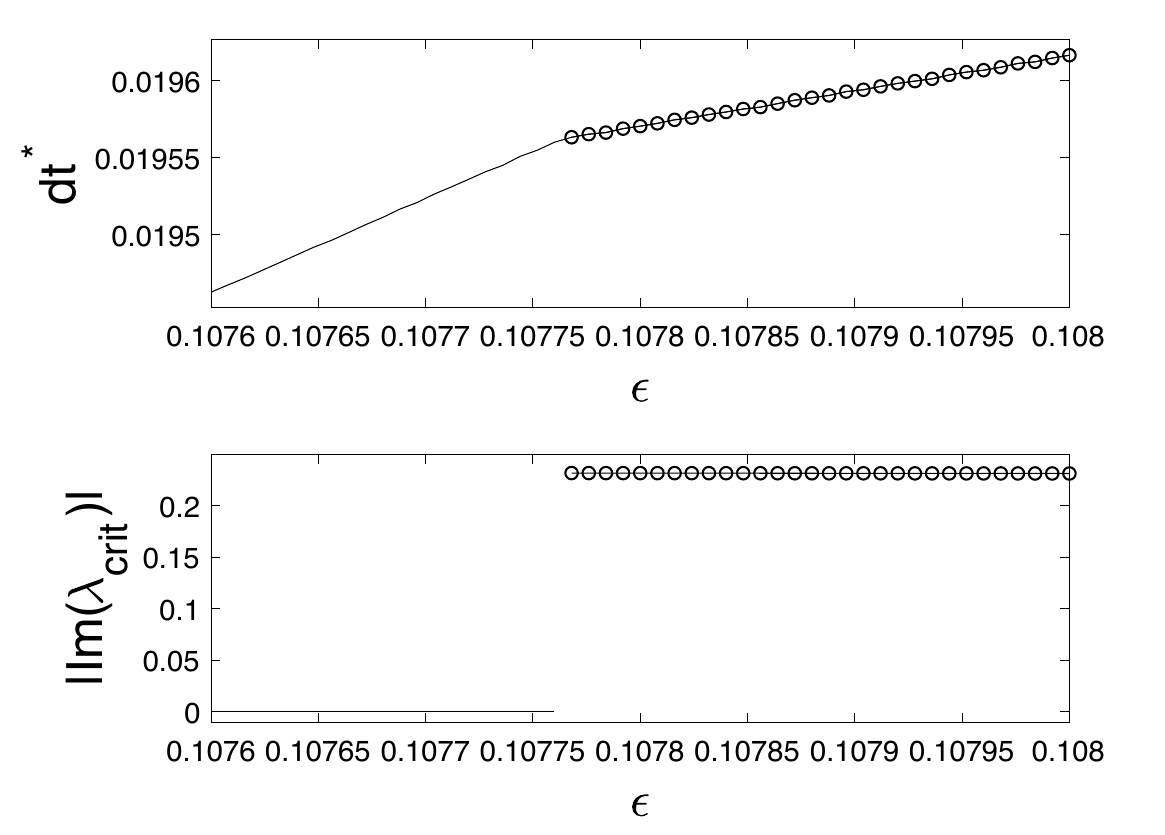}
\includegraphics[width=0.47\linewidth,height=0.445\linewidth]{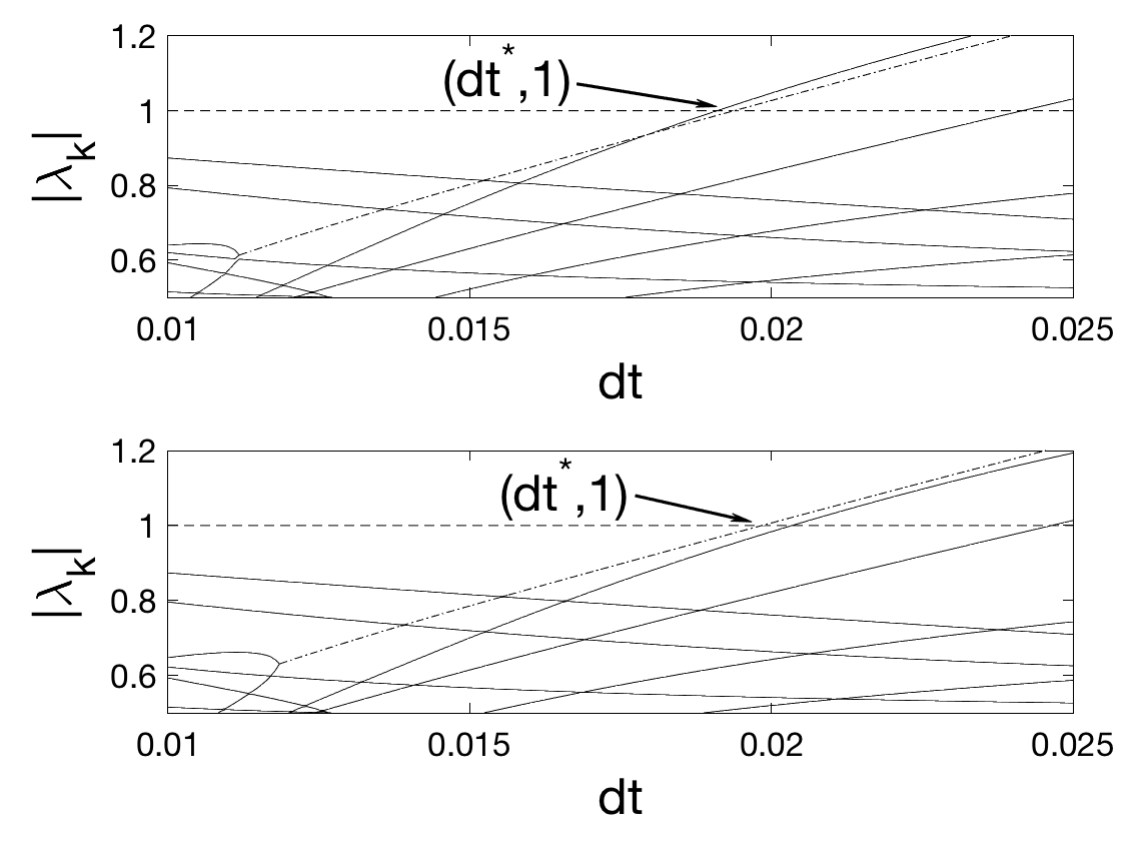}
\caption{PNP-FBV system \eqref{concentration_nondim}--\eqref{phi_bc_nondim_R}
with physical and numerical parameters identical to those used for Figure \ref{Full_Model_tmax_effect} except
for $\epsilon$, which varies.
\underline{Left plot, top:} $dt^*$ versus $\epsilon$ where $dt^*$ is found from the linear stability
analysis.  
There is a corner in the graph at approximately $\epsilon = .107764$.
\underline{Left plot, bottom:}
The magnitude of the imaginary part of the eigenvalue(s) 
on the unit circle.  There is a jump at approximately $\epsilon =  0.107764$.
\underline{Right plot, top:}
The magnitude of the eigenvalues of the
linearized problem versus $dt$ versus for $\epsilon = .107$.
\underline{Right plot, bottom:}
The magnitude of the eigenvalues of the
linearized problem versus $dt$ versus for $\epsilon = .109$.
}
\label{zoom_in_1}
\end{figure}

Figure \ref{zoom_in_1} addresses the cause of the corner
in the graph of $dt^*(\epsilon)$ at $\epsilon \approx .107764$.
The top left
figure presents a closer view of $dt^*$ versus $\epsilon$; there is clearly a corner in the graph.  The bottom left
figure presents the magnitude of the imaginary part of the eigenvalue(s)
that cross the
unit circle.  For $\epsilon$ close to, but smaller than, $\epsilon = .107764$, the instability arises when a single
real-valued eigenvalue crosses the unit circle through the point $-1$.  For $\epsilon$ close to, but larger than, 
$\epsilon = .107764$, the instability arises when a complex conjugate pair of eigenvalues with nonzero
imaginary part cross the unit circle.
The
plots to the right present the magnitudes of the eigenvalues as a function of $dt$ for two values 
of $\epsilon$ close to $\epsilon= .107764 $.  The top plot is for an $\epsilon$ that is close to, but smaller
than, $\epsilon= .107764 $ and the bottom plot is for an $\epsilon$ value that is slightly larger than this
critical value of $\epsilon$.
In both plots, there is a branch which denoted with a dot-dash line.  
This branch corresponds to a pair of complex conjugate eigenvalues; 
following this branch leftward and downward in the figure, one sees that it arose from the collision
of two real-valued eigenvalues (there is a triple junction).  
In the top ($\epsilon = .107$) plot, the branch is to the right of the
branch with one real eigenvalue: the complex pair of eigenvalues are not the cause of the stability restriction.  In the bottom ($\epsilon = .109$)
plot, the two branches have exchanged positions.
If one views a sequence\footnote{The matlab source code is available at \url{https://github.com/daveboat/vssimex_pnp}.  At that site, a curious reader can find a movie of plots like the left plot in Figure \ref{Full_Model_eps_0p12}
as $\epsilon$ increases: {\tt stability\_roots.avi} and {\tt stability\_roots.mov} .} of such plots as $\epsilon$ increases from
$.107$ to $.109$, one sees that both branches are moving rightward but that the branch that
carries the single real eigenvalue is moving rightward at a slightly faster speed; as a result it
overtakes the branch that carries the complex pair of eigenvalues.  The difference
in speeds is the cause of the corner in the graph of $dt^*(\epsilon)$.

\begin{figure}[htb!]
\centering
\includegraphics[width=0.47\linewidth,height=0.458\linewidth]{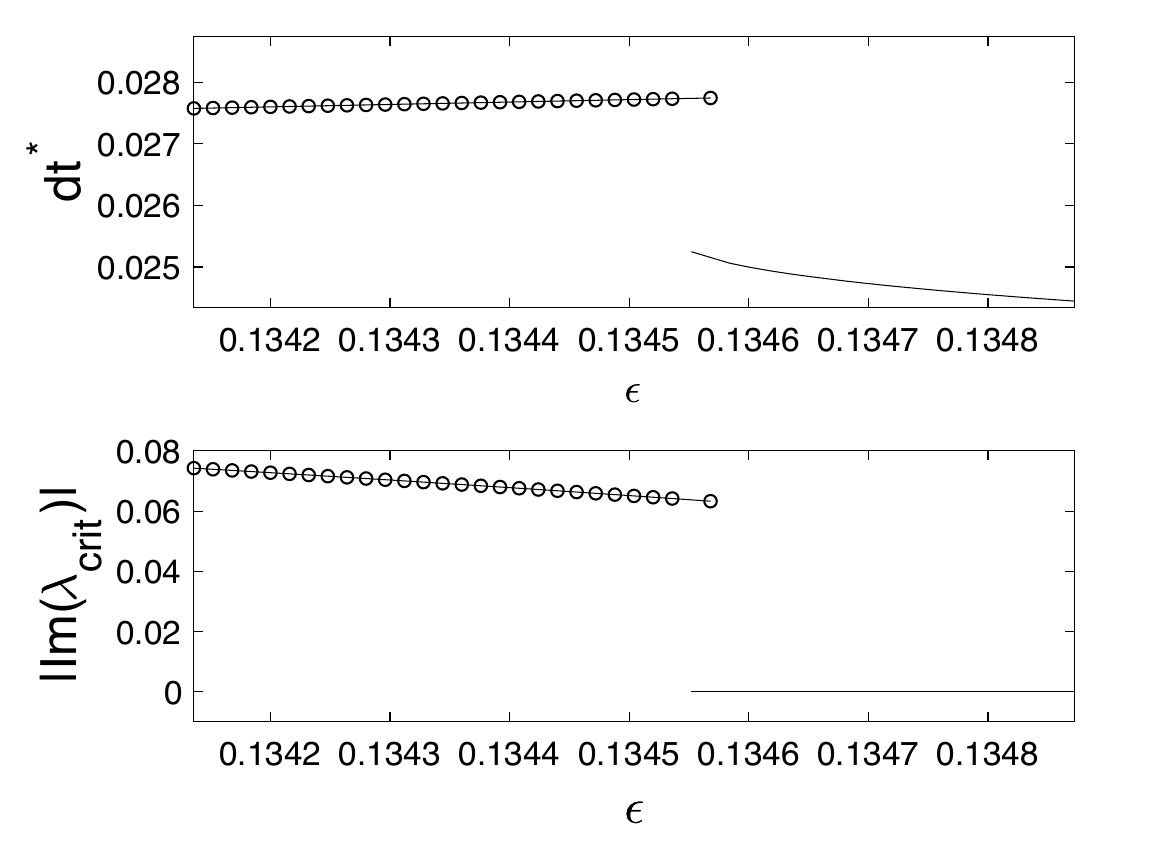}
\includegraphics[width=0.47\linewidth,height=0.455\linewidth]{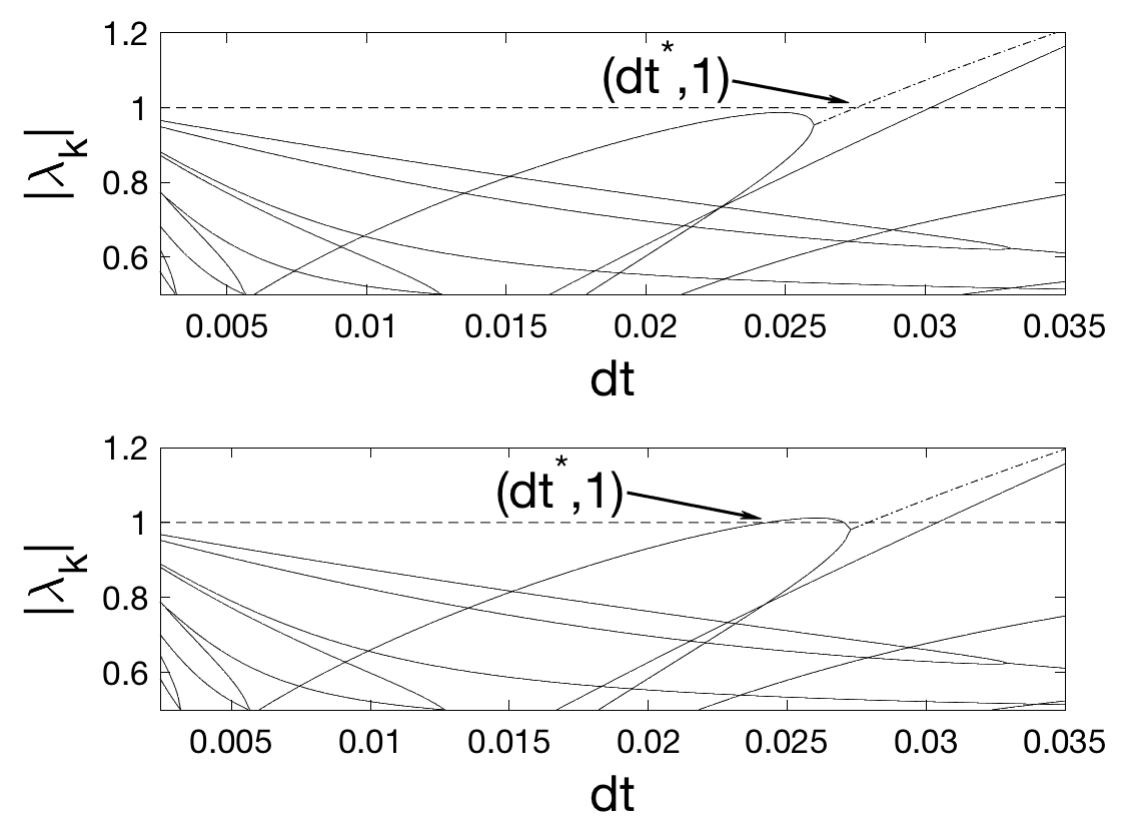}
\caption{PNP-FBV system \eqref{concentration_nondim}--\eqref{phi_bc_nondim_R} 
with physical and numerical parameters identical to those used for Figure \ref{Full_Model_tmax_effect} except
for $\epsilon$, which varies.
\underline{Left plot, top:} $dt^*$ versus $\epsilon$ where $dt^*$ is found from the linear stability
analysis.  
There is a jump at approximately $\epsilon \approx .134504$.
\underline{Left plot, bottom:}
The magnitude of the imaginary part of the eigenvalue(s) on the unit circle.  
\underline{Right plot, top:}
The magnitude of the eigenvalues of the
linearized problem versus $dt$ versus for $\epsilon = .134$.
\underline{Right plot, bottom:}
The magnitude of the eigenvalues of the
linearized problem versus $dt$ versus for $\epsilon = .135$.
}
\label{zoom_in_2}
\end{figure}

Figure \ref{zoom_in_2} is the analogue of Figure \ref{zoom_in_1};
it addresses the cause of the jump
in the graph of $dt^*(\epsilon)$ at the critical value $\epsilon \approx 0.134504$.
From the figure in the left, we see that there is a jump in the stability restriction $dt^*$ and that the
eigenvalues switch from a complex conjugate pair to a single real eigenvalue as $\epsilon$ increases through
the critical value.  The upper right plot presents the magnitude of the 
eigenvalues for a value of $\epsilon$ that is slightly smaller than the critical value and 
the bottom right plot presents them for a value that is slightly larger.   In both plots, we see that the 
triple point, where the branch carrying the complex pair of eigenvalues emerges from the intersection of two branches
carrying single real eigenvalues, is close to the dashed line at height 1.  In the upper plot, we see that the upper branch
(before the triple point) is below the dashed line --- the first eigenvalues to cross the unit circle are the complex pair, for
a larger value of $dt^*$.  However, as $\epsilon$ increases, this upper branch (before the triple point) moves upwards
and it reaches the dashed line when $\epsilon  \approx 0.134504$; at this value of $\epsilon$, the stability restriction $dt^*$
jumps downwards.  After this critical value of $\epsilon$, the stability
restriction is due to a single real eigenvalue crossing the unit circle at $-1$.  

\subsection{Dependence of stability domain on imposed voltage}
\label{v_dpendence}

We next study how the stability restriction, $dt^*$, depends on the
imposed voltage.  
For this, we imposed constant voltages, with values
ranging between $0$ and $3$.  

\begin{figure}[htb!]
\centering
\includegraphics[width=0.47\linewidth,height=0.475\linewidth]{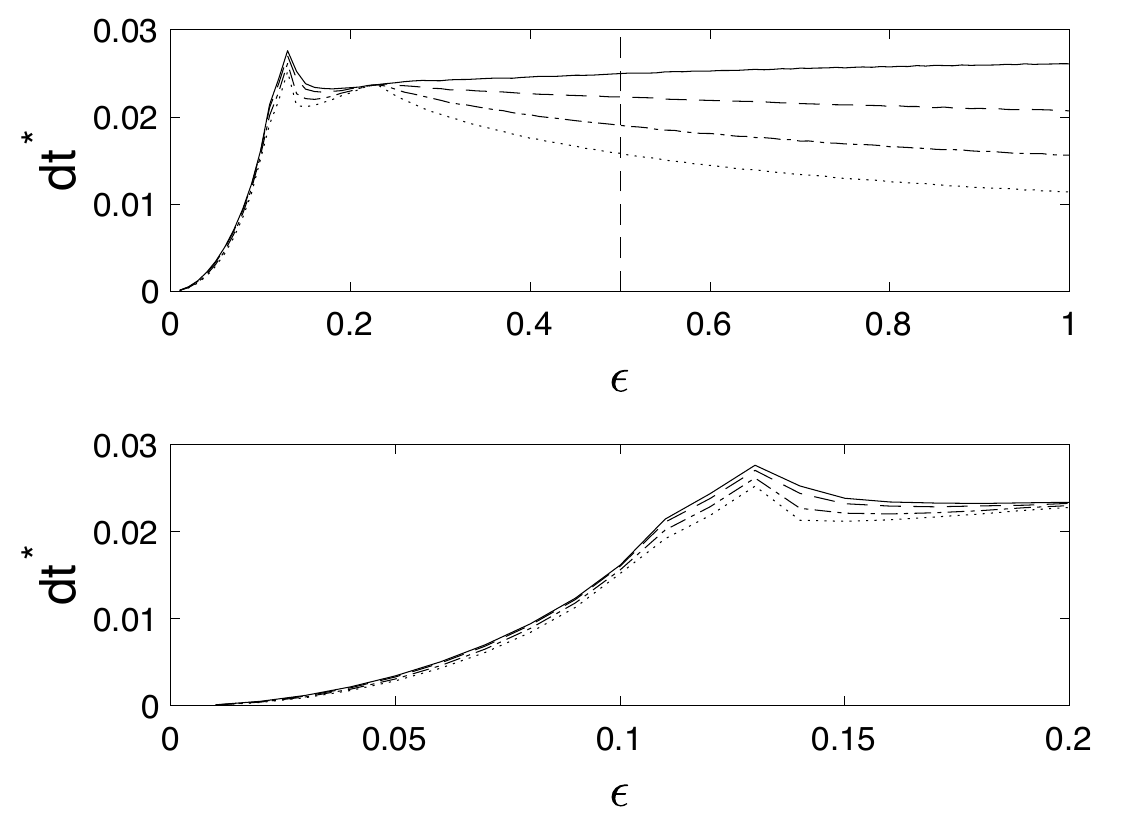}
\includegraphics[width=0.47\linewidth]{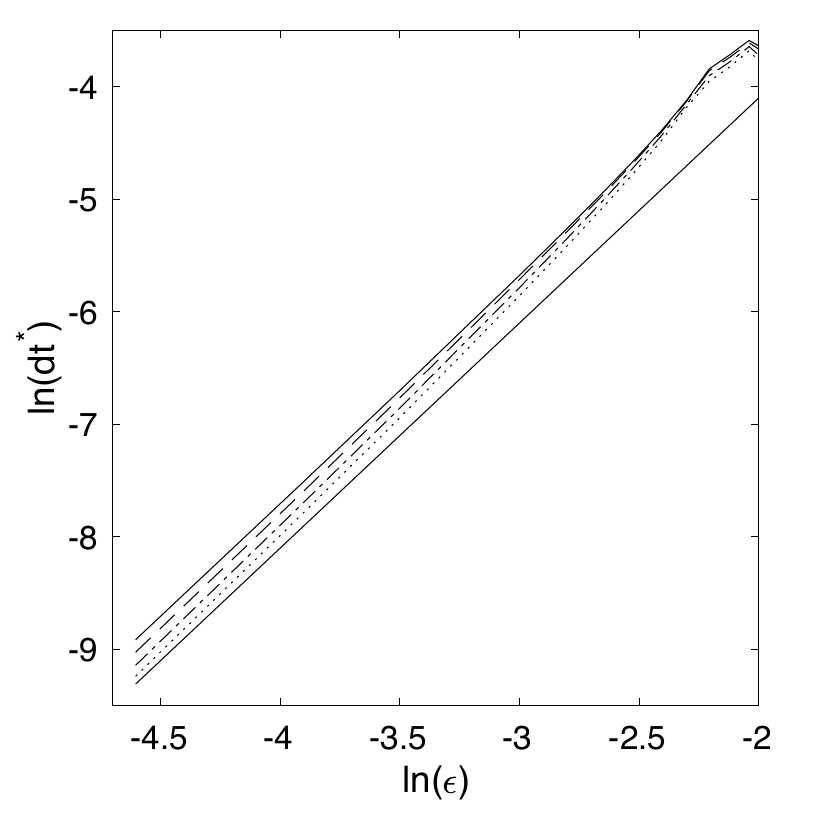}
\caption{PNP-FBV system \eqref{concentration_nondim}--\eqref{phi_bc_nondim_R} 
with physical and numerical parameters identical to those used for Figure \ref{Full_Model_tmax_effect} except
for $\epsilon$, which varies, and the constant imposed  voltage.
Four values of (constant) imposed voltage are considered: $v(t) = 0, 1, 2,$ and $3$.
\underline{Left plots:} The stability restriction $dt^*$ versus $\epsilon$
for $v(t)=0$ (solid), $v(t) = 1$ (dashed), $v(t) = 2$ (dot-dashed),
and $v(t)=3$ (dotted). 
\underline{Right plot:}  Here, $\log(dt^*)$ is plotted versus $\log(\epsilon)$
for the same imposed voltage values and line notations.  The decay is roughly
like $\epsilon^2$; fitting the data yields exponents
2.0677, 2.1145, 2.1540, and 2.1813 for constant imposed voltages $v(t)=0$, $1$,
$2$, and $3$ respectively.  A sight-line corresponding to $dt \propto \epsilon^2$ is provided (solid line).
}
\label{Full_model_dtthresh_vs_epsilon_vary_volt}
\end{figure}

The left plot of Figure
\ref{Full_model_dtthresh_vs_epsilon_vary_volt} presents the stability
domain, $dt^*$, versus $\epsilon$ for the four constant imposed voltages.  We see
that for smaller values of $\epsilon$, $dt^*$ does not appear to be
affected as much by the imposed voltage compared to larger values of
$\epsilon$.  The vertical dashed line in the upper left plot of 
Figure \ref{Full_model_dtthresh_vs_epsilon_vary_volt} indicates
$\epsilon = .5$; it intersects the graphs of $dt^*$ at
the values $0.0250, 0.0223, 0.0191,$ and $0.0158$.
In the companion article \cite{YPD_Part1}, we presented the results of a 
simulation with a time-dependent imposed voltage; the imposed voltage was initially held constant at $0$ and
then transitioned quickly to the value $3$ and held constant.  As demonstrated
in Figure 4 of the companion article \cite{YPD_Part1}, the VSSBDF2 adaptive
time-stepper had its time steps stabilize to $dt \approx 0.0250$ when the
imposed voltage was $0$ and then, after a transient, they stabilized
to $dt \approx 0.0158$ after the imposed voltage was switched to $3$.


The lower left
plot of Figure \ref{Full_model_dtthresh_vs_epsilon_vary_volt} suggests
that $dt^*$ may be proportional to a power of $\epsilon$ for small
values of $\epsilon$.  The right plot of Figure
\ref{Full_model_dtthresh_vs_epsilon_vary_volt} presents $\ln(dt^*)$
versus $\ln(\epsilon)$.  The four graphs appear to be roughly linear
but do not appear to have the same slopes.  All four plots correspond
to $dt^*$ decreasing to zero slightly faster than $\epsilon^2$,
consistent with Table 2 in the companion article \cite{YPD_Part1}.

\subsection{Dependence of stability domain on spatial discretization}
\label{discretization_effect}

\begin{figure}[htb!]
\centering
\includegraphics[width=0.47\linewidth]{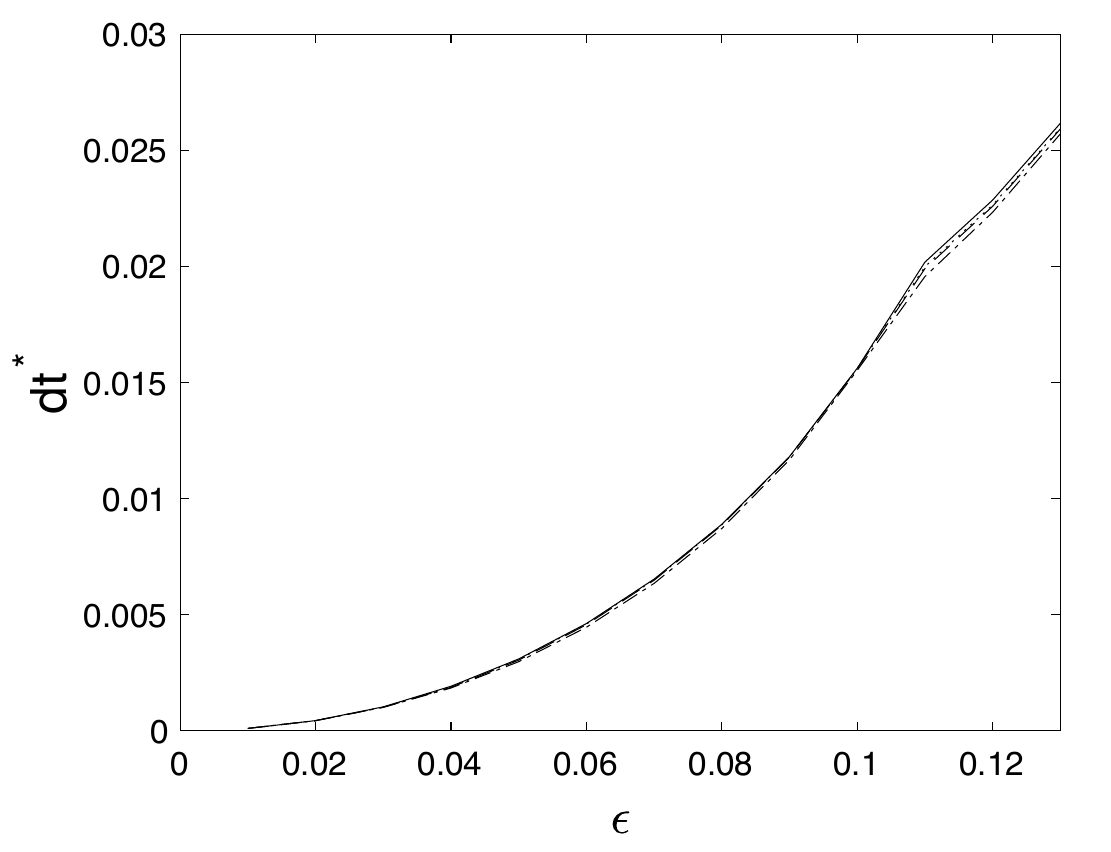}
\includegraphics[width=0.47\linewidth]{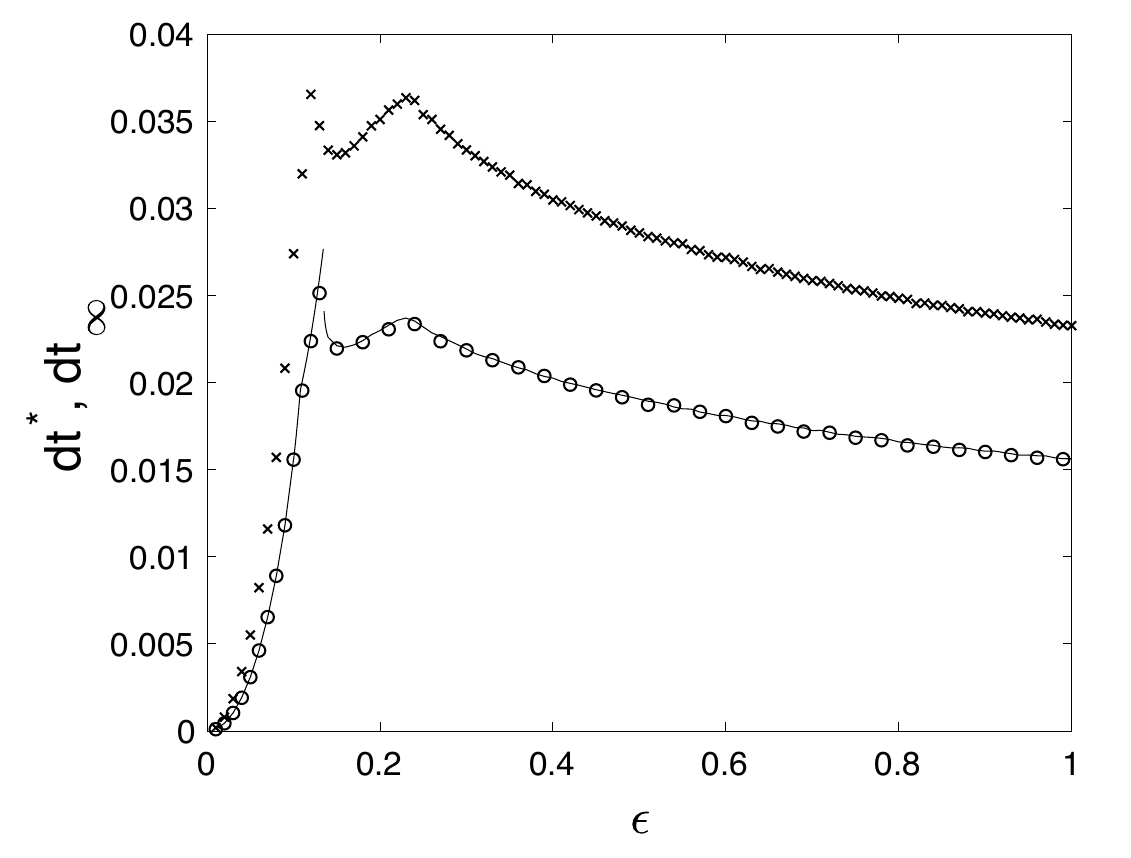}
\caption{PNP-FBV system \eqref{concentration_nondim}--\eqref{phi_bc_nondim_R}
with physical and numerical parameters identical to those used for Figure \ref{Full_Model_tmax_effect} except
for $\epsilon$, which varies.
\underline{Left plot:} For four different
meshes, $dt^*$ is computed and plotted against $\epsilon$.  Solid line: uniform mesh with $dx = 1/90$.
Dot-dashed line: piecewise uniform mesh with $dx = 1/150$ in $[0,1/10]$ and $[9/10,1]$ and $dx = 4/75$ elsewhere.
Dashed line: piecewise uniform mesh with $dx = 1/300$ in $[0,1/10]$ and $[9/10,1]$ and $dx = 2/75$ elsewhere.
Dotted line: piecewise uniform mesh with $dx = 1/450$ in $[0,1/10]$ and $[9/10,1]$ and $dx = 4/225$ elsewhere.
\underline{Right plot:}  The parameters are as in Figure \ref{Full_model_dtthresh_vs_epsilon}.  Solid line: $dt^*$ versus $\epsilon$ where $dt^*$ is found from the linear stability
analysis.  Open circles: $dt_\infty$ as found by the time-step size stabilizing in the VSSBDF2 adaptive
time-stepper with no Richardson extrapolation step.  X-marks: $dt_\infty$ as found by the VSSBDF2 adaptive
time-stepper with a Richardson extrapolation step. 
}
\label{Full_model_dtthresh_vs_epsilon_vary_mesh}
\end{figure}

We next consider the effect of the mesh on the stability restriction.  
The plot in the left of Figure \ref{Full_model_dtthresh_vs_epsilon_vary_mesh} presents stability
restrictions for four different meshes: one
uniform mesh and three piecewise uniform meshes that have finer meshes 
near $x=0, 1$.  
The stability restriction depends on the mesh in a mild manner.  This
is not surprising given that the critical eigenmodes of the linearized scheme
presented in the bottom right plots of Figures \ref{Full_Model_eps_0p05} and \ref{Full_Model_eps_0p12} do not
appear to have structures that need significant spatial resolution. 

We find that the observed dependence is more like that of a reaction--diffusion equation
that that of a diffusion equation.  
Consider
a simple linear PDE: a diffusion equation with a
sink $u_t = D \, u_{xx} - (\frac{1}{\epsilon^2}) \, u$ on $(0,1)$ with 
Dirichlet boundary counditions $u(0,t)=u(1,t)=0$.  Choosing a simple spatial discretization, the resulting system of ODEs is
easily diagonalized.  Using the SBDF2 scheme with the diffusion term handled implicitly and the source term handled
explicitly, one can find the stability restriction $dt^*$ using the
linear stability analysis in Appendix A of \cite{YPD_Part2_2019}.
One finds that that if 
$\epsilon^2 D |\lambda_1| < 3$ then the scheme is conditionally stable
with
$dt^* = 4/(D_1 \lambda_1 + 3/\epsilon^2)$.  Otherwise the scheme is unstable. 
Here $\lambda_1$ is the first negative eigenvalue
of the discretized Laplacian.
For this example, we
see that the stability restriction, $dt^*$, does depend on the mesh because
$\lambda_1$ depends on the mesh.  However, for a fixed $\epsilon$, $dt^*$ does not go to zero as
the mesh is refined: $dt^*$
converges to a positive number.  
If the time-step violates the stability restriction, the fastest growing mode is the eigenvector that
approximates the low frequency eigenfunction $\sin(\pi x)$.   We also note that
$dt^* = 4/(D_1 \lambda_1 + 3/\epsilon^2) \sim \epsilon^2$
for $\epsilon \ll 1$.

\subsection{Effect of Richardson Extrapolation}
\label{RE_effect}
Richardson extrapolation is commonly used to increase the accuracy of
time-stepping, however we have not found much discussion of its possible
effect on numerical stability.    Although we note that in \cite{Eggers_IMEX}, the
authors present an analysis of a stabilized Euler time-stepping scheme
for $u_t = - a \, u$ where $a > 0$.   They demonstrate that their
stabilization parameter $b$ has a threshold value at which their scheme
changes form conditionally stable to unconditionally stable.
They also demonstrate how Richardson
extrapolation affects that threshold value.

The adaptive time-stepper
computes two approximations of the solution at time $t^{n+1}$.
The solution at time $t^{n+1}$ is then taken to equal the ``coarse'' approximation:
$\u^{n+1}_c$; see Section \ref{crash_course}. 
The linear
stability analysis is built upon the solutions satisfying
\eqref{BDF2_IMEX} and so its predictions only apply to the VSSBDF2
adaptive time-stepper when Richardson extrapolation is not used.  

However, as described in \cite{YanThesis,YPD_Part1} and in 
Appendix \ref{ATS}, one could
use $u^{n+1}_c$ and $u^{n+1}_f$ to construct a more accurate
approximation of $u(t^{n+1})$ via Richardson extrapolation.  We find
that when using the VSSBDF2 adaptive time-stepper with Richardson
extrapolation on the PNP-FBV system, the observed behaviour is like
that when Richardson extrapolation was not used: the time-steps
stabilized at a value $dt_\infty$ \cite{YPD_Part2_2019}.  



The plot in the right of Figure \ref{Full_model_dtthresh_vs_epsilon_vary_mesh} presents $dt^*$ and $dt_\infty$
where $dt_\infty$ is found using two different implementations of the VSSBDF2 adaptive time-stepper.
The open circles denote $dt_\infty$ as found by the VSSBDF2 adaptive 
time-stepper with no Richardson extrapolation.  The 
crosses denote $dt_\infty$ as found by the VSSBDF2 adaptive time-stepper with
Richardson extrapolation.  
Because $dt_\infty$ and $dt^*$ agree closely when Richardson extrapolation is not used, we believe
$dt_\infty$ is a good proxy
for the stability restriction $dt^*$ when Richardson extrapolation is used.
In the plot, the crosses are above the circles,
sometimes markedly so,
and for this reason the simulations finish more quickly when Richardson extrapolation is used.
It is also interesting that the plot of $dt_\infty$ versus $\epsilon$ when Richardson extrapolation is used (crosses) has a very similar shape
to the plot of the data when it is not used (circles).  

For the PNP-FBV system, we find that using Richardson extrapolation as part of the adaptive time-stepper leads to 
greater stability.  However, this is problem dependent.
For example, in \cite{YPD_Part2_2019}, 
we find that if one repeats this experiment for the diffusion equation
$u_t = D_1 u_{xx} + D_2 u_{xx}$
then
using Richardson extrapolation
in the VSSBDF2 adaptive time-stepper can lead to {\it less} stability.  
For some choices of $D_1$ and $D_2$,  
 $dt_\infty$ is
smaller 
when 
Richardson extrapolation is used in the VSSBDF2 adaptive time-stepper.
Also, there are choices of $D_1$ and $D_2$ for which the 
SBDF2 time-stepper is unconditionally stable (and so the VSSBDF2 adaptive time-stepper with no
Richardson extrapolation has no stabilization to $dt_\infty$)
but if one uses Richardson extrapolation 
in the adaptive time-stepper then
a limiting step size $dt_\infty$ is observed.  This suggests that, for such parameter choices, using Richardson extrapolation changes the underlying
time-stepping scheme from unconditionally stable to conditionally stable.

\section{Conclusions and Future Work}
\label{conclusions}

In this work, we considered the Poisson-Nernst-Planck equations with generalized Frumkin-Butler-Volmer reaction
kinetics at the electrodes.  
When the VSSBDF2 adaptive time-stepper is being used to study scenarios in which the imposed voltage or the imposed
current is (nearly) constant for long periods of time, the time-step sizes stabilize to a limiting
value and the computed solutions ``nearly'' converge to a steady state.  This 
behaviour is understood by linearizing the numerical scheme about the steady state.   The linearized
scheme is found to be conditionally stable, with a stability restriction that agrees with the time-step at
which that the adaptive time-stepper stabilized.  The stability domain's
dependence on the singular perturbation parameter $\epsilon$ is studied numerically and is found
to have a corner and a jump discontinuity.  The eigenfunctions corresponding to the critical eigenvalues
are studied; the conditional stability is not related to a high-frequency instability.  Using a Richardson
extrapolation step in the adaptive time-stepper appears to stabilize the problem somewhat in that
the limiting time-step is larger.  However other systems are presented for which Richardson extrapolation
can destabilize the scheme.

It would be interesting to see if one can modify, or remove, the conditional stability by using information about the
structure of the steady state.  For example, \cite{Seibold_Shirokoff_Zhou} created an unconditionally stable scheme for a
nonlinear diffusion equation by using bounds on the solution, although the instability being controlled was
due to high frequencies.  

Our methods are not restricted to the PNP-FBV system or to the
VSSBDF2 adaptive time-stepper.  If one is using linear multi-step method 
to study a system that has asymptotically stable
steady states, our approach is relevant.  We expect that it would generalize in a natural manner
to Runge-Kutta methods as well.  A natural next step would be to study the
stability properties of semi-implicit schemes beyond steady states by considering problems that have orbitally
stable special solutions, such as travelling waves, or by considering problems that have
asymptotically stable special solutions, such as self-similar solutions.

\appendix

\section{Overview of adaptive time-stepping scheme} \label{ATS}
Algorithm \ref{adaptive_algorithm} shows our adaptive time-stepping scheme.
This type of error control strategy is discussed in Chapter II.4 of Hairer, Norsett and Wanner \cite{Hairer2009}.

\begin{algorithm}
	\caption{Adaptive time-stepping scheme for a single time step} \label{adaptive_algorithm}
	\begin{algorithmic}
		\State $i \gets 0$ \Comment{Reset loop counter for this time step}
		\State $dt_\text{now} \gets dt_\text{old}$ \Comment{Initial guess at $dt$ for this time step}
		\State $u_c^{n+1} \gets$ TimeStep($dt_\text{now}$) \Comment{Coarse step, TimeStep() using Eq. \eqref{vssbdf2} or \eqref{one_step_IMEX}}
		\State $u_f^{n+1} \gets$ TimeStep($dt_\text{now}/2$) \Comment{Fine step}
		\State $\epsilon_c^i \gets$ Error($dt_\text{now}$, $dt_\text{old}$, $u_c^{n+1}$, $u_f^{n+1}$) \Comment{Error() from equation \eqref{final_error}}
		\While{abs($\epsilon_c^i$ - tol) $>$ range} \Comment{Loop until the error is acceptable}
			\If{$i \geq i_\text{max}$} \Comment{Enforce maximum iterations}
			\State $u_c^{n+1} \gets$ TimeStep($dt_\text{min}$)
			\State $u_f^{n+1} \gets$ TimeStep($dt_\text{min}/2$)
			\State \textbf{break}
			\EndIf
			\State $dt_\text{now} \gets \text{min}\left(\text{max}\left(\left(\frac{\text{tol}}{\epsilon_c^i}\right)^{1/p}, \eta_\text{min}\right), \, \eta_\text{max}\right)dt_\text{now}$ \Comment{Update $dt$}
			\If{$dt_\text{now} > dt_\text{max}$} \Comment{Enforce maximum time step}
				\State $u_c^{n+1} \gets$ TimeStep($dt_\text{max}$)
				\State $u_f^{n+1} \gets$ TimeStep($dt_\text{max}/2$)
				\State \textbf{break}
			\EndIf
			\State $u_c^{n+1} \gets$ TimeStep($dt_\text{now}$)
			\State $u_f^{n+1} \gets$ TimeStep($dt_\text{now}/2$)
			\State $\epsilon_c^{i+1} \gets$ Error($dt_\text{now}$, $dt_\text{old}$, $u_c^{n+1}$, $u_f^{n+1}$) \Comment{Update error estimate}
			\State $i \gets i+1$
		\EndWhile
		\If{Richardson extrapolation is used}
		\State $u^{n+1} \gets \alpha u_c^{n+1} + \beta u_f^{n+1}$ \Comment{$\alpha$ and $\beta$ are defined in equation \eqref{final_extrapolation}}
		\Else
		\State $u^{n+1} \gets u_c^{n+1}$
		\EndIf
	\end{algorithmic}
\end{algorithm}
The $p$ in the time-step update formula in the algorithm is the
order of the local truncation error; $p=3$ for VSSBDF2.  
Unless noted otherwise, the simulations in this article used
$tol = 10^{-6}$, $range = tol/3$, $\eta_\text{max} = 1.1$, $\eta_\text{min} = .9$,
$dt_\text{max} = 1$ and $dt_\text{min} = 10^{-8}$.

Richardson extrapolation uses a linear combination of $u_f^{n+1}$ and $u_c^{n+1}$ to construct an improved approximation $u^{n+1}$ which has a smaller truncation error.  Specifically,
$
u^{n+1}=\alpha u^{n+1}_c + \beta u^{n+1}_f
$
with coefficients
\begin{equation}
\label{final_extrapolation}
\alpha=-\frac{dt_\text{old}+3 \, dt_\text{now}}{7 \, dt_\text{old}+5 \, dt_\text{now}}
\quad \mbox{and} \quad
\beta=8\, \frac{dt_\text{old}+dt_\text{now}}{7 \, dt_\text{old}+5 \, dt_\text{now}}.
\end{equation}
The local truncation error for $u^{n+1}$ is one order higher than the local truncation errors for $u^{n+1}_c$ and $u^{n+1}_f$ \cite{YanThesis,YPD_Part1}. Note that if $dt_\text{now}=dt_\text{old}$, then \eqref{final_extrapolation} reduces to the standard Richardson extrapolation formula for second-order schemes.

Since we are using a two-step time-stepping scheme, for the first time-step, we use a one-step semi-implicit scheme 
\begin{equation} \label{one_step_IMEX}
\frac{1}{dt}\left(u^{1} - u^{0} \right)
= f(u^0) + g(u^{1}),
\end{equation}
along with the error 
estimate $\epsilon^1_c = (4/3) (u^{1}_c-u^{1}_f )$, the time-step update formula
with $p=2$, and the extrapolation formula $u^{1}=2u^1_f-u^1_c$.

\section{Acknowledgements} 

Research supported in part by NSERC grant OGP06617.

We thank Greg Lewis, Keith Promislow, Steve Ruuth, Adam Stinchcombe, and Brian Wetton
for helpful conversations and encouragement.  We thank the reviewers
of an earlier version of this article for their careful, thorough
comments and suggestions.



\end{document}